# HIGH BREAKDOWN POINT ROBUST REGRESSION WITH CENSORED DATA

By Matías Salibian-Barrera[1] and Víctor J. Yohai[2]

*University of British Columbia and Universidad de Buenos Aires*

In this paper, we propose a class of high breakdown point estimators for the linear regression model when the response variable contains censored observations. These estimators are robust against high-leverage outliers and they generalize the LMS (least median of squares), S, MM and $\tau$-estimators for linear regression. An important contribution of this paper is that we can define consistent estimators using a bounded loss function (or equivalently, a redescending score function). Since the calculation of these estimators can be computationally costly, we propose an efficient algorithm to compute them. We illustrate their use on an example and present simulation studies that show that these estimators also have good finite sample properties.

**1. Introduction.** Consider the linear regression model

$$(1.1) \qquad y_i = \boldsymbol{\beta}_0' \mathbf{x}_i + u_i, \qquad i = 1, \ldots, n,$$

where $u_i$ are i.i.d. errors, and the covariates $\mathbf{x}_i \in \mathbb{R}^p$ are independent from the errors. When there is an intercept the first component of $\mathbf{x}_i$ is set to 1. In this paper, we study the problem of robust estimation of $\boldsymbol{\beta}_0$ when the response variable is censored. Miller [12] studied least squares estimators (LS) for censored responses. He proposed to modify the classical LS estimator

$$(1.2) \qquad \hat{\boldsymbol{\beta}}_n = \underset{\boldsymbol{\beta} \in \mathbb{R}^p}{\arg\min} \sum_{i=1}^{n}(y_i - \boldsymbol{\beta}' \mathbf{x}_i)^2 = \underset{\boldsymbol{\beta} \in \mathbb{R}^p}{\arg\min} E_{F_{n\boldsymbol{\beta}}}[u^2],$$

Received November 2006; revised April 2007.
[1]Supported by a Natural Sciences and Engineering Research Council of Canada Discovery grant.
[2]Supported in part by Grant X-094 from the Universidad de Buenos Aires, Grant PIP 5505 from CONICET, Argentina and Grant PICT 21407 from ANPCYT, Argentina.

replacing the empirical distribution of the residuals $F_{n\boldsymbol{\beta}}$ with the corresponding Kaplan–Meier (KM) estimator $F^*_{n\boldsymbol{\beta}}$ (Kaplan and Meier [9]). Unfortunately, the resulting estimator is not consistent in general and the iterative algorithm to compute it may have several or no solutions.

Buckley and James [2] studied a different extension of LS to censored response variables by modifying the LS scores equations

$$(1.3) \qquad \sum_{i=1}^{n}(y_i - \hat{\boldsymbol{\beta}}'_n \mathbf{x}_i)\mathbf{x}_i = 0,$$

using a conditional distribution approach. This proposal replaces censored residuals by their estimated conditional expectation given that the response is larger than the recorded (censored) value. The conditional expectation is estimated using $F^*_{n\hat{\boldsymbol{\beta}}_n}$. James and Smith [7] and Lai and Ying [10] showed that this estimator is consistent.

A different approach is proposed by Stute [19, 20] and Sellero, Manteiga and Van Keilegom [18]. They propose to apply Kaplan–Meyer to the responses instead to the residuals. The shortcomings of this approach is that they require stronger assumptions on the censoring variable and that the proposed estimates are not regression equivariant.

In recent years there has been some interest in extending robust regression estimators to the case of censored response variables. Ritov [13] studied a generalization of Bukley and James' proposal for robust estimators. He considered monotone nondecreasing score functions $\psi$ (that correspond to unbounded loss functions $\rho$) and showed that under certain regularity conditions there exists a sequence of $\sqrt{n}$-consistent solutions to the estimating equations. This sequence is also asymptotically normal. Unfortunately, since these estimators are based on an unbounded loss function $\rho$ they are not robust against high-leverage outliers. More recently, Lai and Ying [11] extended the conditional expectation approach of Bukley and James to M-regression estimators for censored and truncated data. Their proposal also requires a monotone score function.

If we allow for a redescending score function $\psi$ (equivalently, a bounded loss function $\rho$), then the estimating equations may have several solutions with different robustness properties. Moreover, if we define a robust estimator as the solution to a minimization problem similar to (1.2) but replacing the squared residuals with $\rho(u)$ for a bounded loss function $\rho$, then this estimator may not be consistent (Lai and Ying [10, 11]). Hence, unlike in the uncensored regression model, we do not have a way to identify which solutions of the redescending score equations are not affected by the outliers.

In this paper, we extend the approach of Bukley and James and Ritov to M-estimators with bounded loss functions $\rho$. We achieve this by proposing an estimator that is the solution to a minimization problem that has a consistent and robust solution. In particular, we obtain extensions of the LMS

ROBUST CENSORED REGRESSION 3(see Rousseeuw [14]), S (see Rousseeuw and Yohai [16]), MM-estimators (see Yohai [22]) and $\tau$-estimators (see Zamar and Yohai [23]). We show that these estimators are Fisher and $\sqrt{n}$-consistent, asymptotically normal, and that they have high breakdown point.

It is important to realize that when there are censored observations the breakdown point of an estimator maybe much lower than in the uncensored case. For example, in the location model the worst contamination occurs when all the censored observations are between the outliers and the "good" noncensored points. Suppose that we have a fraction $\epsilon$ of outliers going to $+\infty$ and a proportion $\lambda$ of censored observations. Since the KM estimator distributes the mass of the censored observations among the noncensored points to their right (Efron [4]), in this case the mass given to the outliers by the KM estimators will be $\gamma = \lambda + \epsilon$. Consequently, the sample median will not break if $\gamma < 1/2$, or equivalently, if $\epsilon < 1/2 - \lambda = \eta$. It follows that the breakdown point of the median is equal to $\eta$, which is less than $1/2$ when there are censored observations.

The rest of this paper is organized as follows. Section 2 contains our main definitions. The robustness properties of our proposal are discussed in Section 3 and their asymptotic properties in Section 4. In Section 5, we present an algorithm to compute these estimators. An example with real-life data is given in Section 6 and the results of a Monte Carlo experiment are discussed in Section 7. The proofs of the theorems are given in the Appendix while those for the lemmas can be found in a technical report by Salibian-Barrera and Yohai [17].

**2. Robust estimators.** Consider the linear regression model (1.1). We assume that the sample may be right-censored, that is, there are unobservable random variables $c_1, \ldots, c_n$ independent from the errors $u_i$'s such that we observe $y_i^* = \min(y_i, c_i)$ for $i = 1, \ldots, n$. In other words, the observed data is $\mathbf{z}_i = (y_i^*, \mathbf{x}_i', \delta_i)'$, $i = 1, \ldots, n$, where $\delta_i = I\{y_i \leq c_i\}$, and $I\{A\}$ is the indicator function of the event $A$.

When the scale of the residuals is known, regression M-estimators for uncensored observations are defined by

$$(2.1) \qquad \hat{\boldsymbol{\beta}}_n = \arg\min_{\boldsymbol{\beta} \in \mathbb{R}^p} \frac{1}{n} \sum_{i=1}^n \rho(r_i(\boldsymbol{\beta})) = \arg\min_{\boldsymbol{\beta} \in \mathbb{R}^p} E_{F_{n\boldsymbol{\beta}}}[\rho(u)],$$

where $F_{n\boldsymbol{\beta}}$ is the empirical distribution of the residuals $r_i(\boldsymbol{\beta}) = y_i - \boldsymbol{\beta}'\mathbf{x}_i$, and $\rho: \mathbb{R} \to \mathbb{R}^+$ is a function satisfying:

P1. $\rho(0) = 0$ and $\rho$ is continuous at 0.
P2. $\rho(-u) = \rho(u)$ for $u > 0$.
P3. $\rho$ is monotone nondecreasing on $u > 0$.
P4. $\sup_u \rho(u) = a < +\infty$.

(See Huber [6].) If $\psi(u) = \partial \rho(u)/\partial u$ then the estimator $\hat{\boldsymbol{\beta}}_n$ also satisfies the following vector equation:

$$(2.2) \qquad \frac{1}{n}\sum_{i=1}^{n} \psi(r_i(\boldsymbol{\beta}))\mathbf{x}_i = E_{H_{n\boldsymbol{\beta}}}[\psi(u)\mathbf{x}] = \mathbf{0},$$

where $H_{n\boldsymbol{\beta}}$ is the empirical distribution of the vectors $(r_i(\boldsymbol{\beta}), \mathbf{x}_i')' \in \mathbb{R}^{p+1}$, $i = 1, \ldots, n$.

Since not all the residuals $r_i(\boldsymbol{\beta})$ are observed in the presence of censoring, we can define the censored residuals by $r_i^*(\boldsymbol{\beta}) = y_i^* - \boldsymbol{\beta}'\mathbf{x}_i$. Note that $r_i^*(\boldsymbol{\beta}) = \min(r_i(\boldsymbol{\beta}), c_i - \boldsymbol{\beta}'\mathbf{x}_i)$, and therefore we can think of the $r_i^*(\boldsymbol{\beta})$ as censored observations of $r_i(\boldsymbol{\beta})$ with censoring variables $c_i - \boldsymbol{\beta}'\mathbf{x}_i$, $i = 1, \ldots, n$. Then, in the case of a censored response variable, one way to generalize (2.1) is to replace it by

$$(2.3)\quad \hat{\boldsymbol{\beta}}_n = \arg\min_{\boldsymbol{\beta}\in\mathbb{R}^p} \frac{1}{n}\sum_{i=1}^{n} E[\rho(r_i(\boldsymbol{\beta}))|\mathbf{z}_i] = \arg\min_{\boldsymbol{\beta}\in\mathbb{R}^p} \frac{1}{n}\sum_{i=1}^{n} E_{F_{\boldsymbol{\beta}}}[\rho(u)|\mathbf{w}_i(\boldsymbol{\beta})],$$

where $F_{\boldsymbol{\beta}}$ is the distribution of the residuals $r(\boldsymbol{\beta})$, $\mathbf{w}_i(\boldsymbol{\beta}) = (r_i^*(\boldsymbol{\beta}), \delta_i)$ and

$$E_{F_{\boldsymbol{\beta}}}(\rho(u)|\mathbf{w}_i(\boldsymbol{\beta})) = \begin{cases} \rho(r_i^*(\boldsymbol{\beta})), & \text{if } \delta_i = 1, \\ \int_{r_i^*(\boldsymbol{\beta})}^{\infty} \rho(u)\,dF_{\boldsymbol{\beta}}(u)/[1 - F_{\boldsymbol{\beta}}(r_i^*(\boldsymbol{\beta}))], & \text{if } \delta_i = 0. \end{cases}$$

Intuitively, to obtain (2.3) from (2.1), for each censored observation we replace the term $\rho(r_i(\boldsymbol{\beta}))$ in (2.1) by the conditional expectation of $\rho(u)$ given that the (actual but unobserved) residual is larger than or equal to the observed censored residual $r_i^*(\boldsymbol{\beta})$.

The score equations in (2.2) can also be similarly modified to obtain

$$(2.4) \qquad \frac{1}{n}\sum_{i=1}^{n} E_{F_{\boldsymbol{\beta}}}[\psi(u)|\mathbf{w}_i(\boldsymbol{\beta})]\mathbf{x}_i = \mathbf{0}.$$

Since the distribution of the residuals $F_{\boldsymbol{\beta}}$ in (2.3) and (2.4) is unknown, we can estimate it with the Kaplan–Meier estimator $F_{n\boldsymbol{\beta}}^*$ based on $r_i^*(\boldsymbol{\beta})$.

To guarantee consistency of the estimator defined by

$$(2.5) \qquad \hat{\boldsymbol{\beta}}_n = \arg\min_{\boldsymbol{\beta}\in\mathbb{R}^p} \frac{1}{n}\sum_{i=1}^{n} E_{F_{n\boldsymbol{\beta}}^*}[\rho(u)|\mathbf{w}_i(\boldsymbol{\beta})]$$

we need that $F_{n\boldsymbol{\beta}}^*$ be consistent to $F_{\boldsymbol{\beta}}$ for all $\boldsymbol{\beta}$. Let $F$ and $D$ be the distribution functions of the errors $u_i$ and censoring variables $c_i$, $i = 1, \ldots, n$, respectively. Let $\tau_F = \inf\{u : F(u) = 1\}$ and let $\tau_D$ be defined similarly. In what follows we will assume that:

R1. $\tau_F < \tau_D$, or $\tau_F = \tau_D = \infty$, or $\tau_F = \tau_D$ and $\tau_F$ is a continuity point of $F$.

R2. $F$ and $D$ do not have jumps in common.

Under these conditions, a sufficient condition for the KM estimator to be consistent is the independence between the uncensored variables and the censoring times (see, e.g., Breslow and Crowley [1]). When $\boldsymbol{\beta} = \boldsymbol{\beta}_0$ we have $r_i(\boldsymbol{\beta}_0) = u_i$ which are independent from the corresponding censoring times $c_i - \boldsymbol{\beta}_0' \mathbf{x}_i$ because we have assumed that the errors are independent from the $c_i$'s and the $\mathbf{x}_i$'s. However, for $\boldsymbol{\beta} \neq \boldsymbol{\beta}_0$ it is not generally true that $r_i(\boldsymbol{\beta})$ is independent from $c_i - \boldsymbol{\beta}' \mathbf{x}_i$, $i = 1, \ldots, n$. Hence, we can only guarantee the consistency of $F_{n\boldsymbol{\beta}}^*$ to $F_{\boldsymbol{\beta}}$ when $\boldsymbol{\beta} = \boldsymbol{\beta}_0$. Therefore, the estimator defined in (2.5) may not be consistent (Lai and Ying [10, 11]).

On the other hand, note that the estimator $\hat{\boldsymbol{\beta}}_n$ defined as the solution to

$$(2.6) \qquad \frac{1}{n} \sum_{i=1}^{n} E_{F_{n\boldsymbol{\beta}}^*}[\psi(u) | \mathbf{w}_i(\boldsymbol{\beta})] \mathbf{x}_i = \mathbf{0},$$

is Fisher consistent. In fact, $F_{n\boldsymbol{\beta}_0}^* \to F_{\boldsymbol{\beta}_0}$ and therefore

$$\frac{1}{n} \sum_{i=1}^{n} E_{F_{n\boldsymbol{\beta}_0}^*}[\psi(u) | \mathbf{w}_i(\boldsymbol{\beta}_0)] \mathbf{x}_i \to E_{H_0}(\psi(u)\mathbf{x}) = \mathbf{0},$$

where $H_0$ is the joint distribution of $(u, \mathbf{x}')'$. It is important to note that, unlike in the uncensored regression case, equations (2.5) and (2.6) are not equivalent: we cannot obtain (2.6) by differentiating (2.5) because $F_{n\boldsymbol{\beta}}^*$ depends on $\boldsymbol{\beta}$.

M-estimators defined by (2.6) were first proposed by Ritov [13] and further studied by Lai and Ying [11] when $\psi(u)$ is monotone (which corresponds to a convex $\rho$). However, it is well known that M-estimators with monotone $\psi$ functions are only robust against low leverage outliers. As mentioned in the Introduction, the main difficulty in using a redescending $\psi$ in (2.6) is that in general this equation may have several solutions with different robustness properties. Although in the uncensored regression model this difficulty can be avoided by defining the estimator as the solution to the minimization problem (2.1), the corresponding minimization in the censored case (2.5) does not in general yield a consistent estimator. In other words, (2.5) cannot be used to select a consistent solution of (2.6). For this reason, in the next subsection we will define robust M-estimators as the solution of a minimization problem using a bounded loss function $\rho$ that has a consistent sequence of solutions.

2.1. *Consistent M-estimators.* First note that to obtain scale equivariant regression estimators, we need to standardize the residuals in the estimating equations using a robust error scale estimator $s_n$.

Let $\rho : \mathbb{R} \to \mathbb{R}^+$ satisfy regularity conditions P1–P4 above. For each $\boldsymbol{\beta}$ and $\boldsymbol{\gamma}$ in $\mathbb{R}^p$ define

$$(2.7) \qquad C_n(\boldsymbol{\beta}, \boldsymbol{\gamma}) = \frac{1}{n} \sum_{i=1}^{n} E_{F^*_{n\beta}} \left[ \rho\left(\frac{u - \boldsymbol{\gamma}' \mathbf{x}_i}{s_n}\right) \Big| \mathbf{w}_i(\boldsymbol{\beta}) \right],$$

where $s_n$ is a robust scale estimator of the residuals. For each $\boldsymbol{\beta} \in \mathbb{R}^p$ let

$$(2.8) \qquad \hat{\boldsymbol{\gamma}}_n(\boldsymbol{\beta}) = \arg\min_{\boldsymbol{\gamma} \in \mathbb{R}^p} C_n(\boldsymbol{\beta}, \boldsymbol{\gamma}).$$

Note that $\hat{\boldsymbol{\gamma}}_n(\boldsymbol{\beta})$ can be considered an M-estimator of regression of the residuals $r_i(\boldsymbol{\beta})$ on the covariates $\mathbf{x}_i$. Since $F^*_{n,\boldsymbol{\beta}_0}$ is a consistent estimator of $F_{\boldsymbol{\beta}_0}$, the distribution of the $u_i$'s, and since the errors are independent of the $\mathbf{x}_i$'s, it is reasonable to expect that $\hat{\boldsymbol{\gamma}}_n(\boldsymbol{\beta}_0) \to \mathbf{0}$. This can be formally proved with similar arguments to those used in the proof of Theorem 5 below. Therefore, we define an estimator of $\boldsymbol{\beta}_0$ by the equation

$$(2.9) \qquad \hat{\boldsymbol{\gamma}}_n(\hat{\boldsymbol{\beta}}_n) = \mathbf{0}.$$

To avoid existence problems, we can alternatively define $\hat{\boldsymbol{\beta}}_n$ as

$$(2.10) \qquad \hat{\boldsymbol{\beta}}_n = \arg\min_{\boldsymbol{\beta} \in \mathbb{R}^p} [\hat{\boldsymbol{\gamma}}_n(\boldsymbol{\beta})' \mathbf{A}_n \hat{\boldsymbol{\gamma}}_n(\boldsymbol{\beta})],$$

where $\mathbf{A}_n = \mathbf{A}_n(\mathbf{x}_1, \ldots, \mathbf{x}_n)$ is any robust equivariant estimator of the covariance matrix of the explanatory variables $\mathbf{x}_i$, $1 \le i \le n$. The covariance matrix $\mathbf{A}_n$ is needed to maintain the affine equivariance of the estimator.

As an illustration of the difference between using (2.5) and (2.9) to define a robust estimator, in Figure 1 we plot $\|\hat{\boldsymbol{\gamma}}_n(\boldsymbol{\beta})\|$ and the score equations (2.5) as a function of $\boldsymbol{\beta}$ for a data set of $n = 200$ observations with $\boldsymbol{\beta}_0 = 1.5$ and a probability of censoring of approximately 32%. These data were generated following the same model we used in our simulation study described in Section 7. Note that although the score equation has two distinct solutions and only one is close to the true value of $\boldsymbol{\beta}_0 = 1.5$, our proposed optimization problem has a unique minimum and this minimum is close to $\boldsymbol{\beta}_0$. This definition may be considered an extension of Ritov's M-estimators for censored data to the case of bounded $\rho$ functions. In particular, note that $\hat{\boldsymbol{\beta}}_n$ satisfies equation (2.6) with $\psi(u) = \rho'(u)$. It follows that this estimator will have the same asymptotic properties as the estimators considered in Ritov [13].

2.2. *S-estimators.* The scale estimator $s_n$ in (2.7) may be chosen to be the scale of the residuals of an initial (and scale-equivariant) estimator that does not require a scale estimator itself. One class of estimates that satisfies this is the class of S-estimators (Rousseeuw and Yohai [16]). We can extend this class of estimators to the case of censored observations following the same principle as above, that is, for each $\boldsymbol{\beta}$ we fit an S-estimate to the residuals $r_i^*(\boldsymbol{\beta})$, and find the $\boldsymbol{\beta}$ whose residuals have the "smallest" S-estimator (i.e., the one with the smallest norm).

Let $\rho_1$ satisfy regularity conditions P1–P4 and let $b = E_F[\rho_1(u)]$ where $F$ is the distribution of the errors $u_i$ in (1.1). Define the M-scale $S_n(\boldsymbol{\beta}, \boldsymbol{\gamma})$ by

$$(2.11) \qquad \frac{1}{n} \sum_{i=1}^{n} E_{F_{n\beta}^*}\left[\rho_1\left(\frac{u - \boldsymbol{\gamma}'\mathbf{x}_i}{S_n(\boldsymbol{\beta}, \boldsymbol{\gamma})}\right) \bigg| \mathbf{w}_i(\boldsymbol{\beta})\right] = b$$

and let

$$(2.12) \qquad \hat{\boldsymbol{\gamma}}_n(\boldsymbol{\beta}) = \arg\min_{\boldsymbol{\gamma} \in \mathbb{R}^p} S_n(\boldsymbol{\beta}, \boldsymbol{\gamma}).$$

Note that $\hat{\boldsymbol{\gamma}}_n(\boldsymbol{\beta})$ is the S-estimator of regression of the residuals $(r_i^*(\boldsymbol{\beta}), \mathbf{x}_i')'$, $i = 1, \ldots, n$. We define the S-regression estimator for censored responses as the vector $\tilde{\boldsymbol{\beta}}_n$ such that

$$(2.13) \qquad \hat{\boldsymbol{\gamma}}_n(\tilde{\boldsymbol{\beta}}_n) = \mathbf{0}.$$

As before, to avoid existence problems, the following definition is also natural:

$$\tilde{\boldsymbol{\beta}}_n = \arg\min_{\boldsymbol{\beta} \in \mathbb{R}^p} [\hat{\boldsymbol{\gamma}}_n(\boldsymbol{\beta})' \mathbf{A}_n \hat{\boldsymbol{\gamma}}_n(\boldsymbol{\beta})],$$

where $\mathbf{A}_n = \mathbf{A}_n(\mathbf{x}_1, \ldots, \mathbf{x}_n)$ is any robust equivariant estimator of the covariance matrix of the covariates $\mathbf{x}_i$.

A robust residual scale estimate $s_n$ can be defined by

$$(2.14) \qquad s_n = S_n(\tilde{\boldsymbol{\beta}}_n, \hat{\boldsymbol{\gamma}}_n(\tilde{\boldsymbol{\beta}}_n)).$$

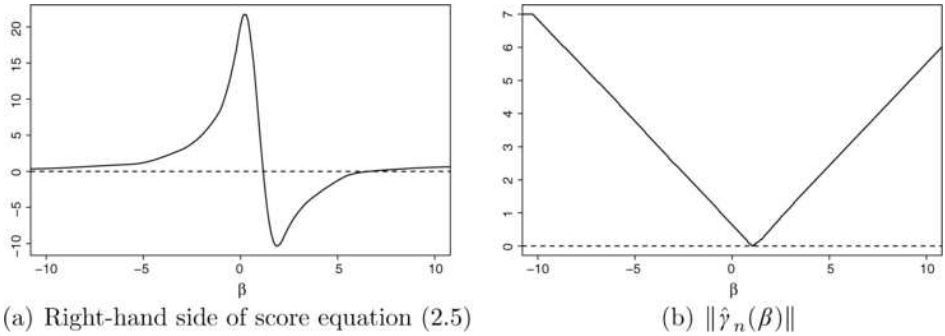

(a) Right-hand side of score equation (2.5)  (b) $\|\hat{\gamma}_n(\beta)\|$

FIG. 1. *Panel (*a*) shows an example where the score equations (2.5) have two roots with only one of them close to $\boldsymbol{\beta}_0 = 1.5$ whereas panel (*b*) shows that, for the same data set, the objective function of (2.10) has a unique minimum close to $\boldsymbol{\beta}_0$.*

In particular, we can obtain a consistent version of the LMS using as $\rho_1$ a jump function

$$\rho_1(u) = \begin{cases} 0, & \text{if } |u| < 1, \\ 1, & \text{if } |u| \geq 1, \end{cases} \tag{2.15}$$

and $b = 1/2$ in equation (2.11) above.

In Section 3, we will show that the choice $b = \sup_u \rho_1(u)/2$ yields regression estimators with high breakdown point. However, we know from the uncensored case that S-estimators cannot combine high breakdown point with high efficiency for normal errors (see Hössjer [5]). To overcome this problem, in the next subsection we will extend to the censored case a class of estimators that can achieve simultaneous high efficiency and high breakdown point.

2.3. *MM-estimators.* Yohai [22] proposed a class of estimators, called MM-estimators, that simultaneously have breakdown point 50% and high efficiency for normal errors. In this subsection, we extend this class of estimators to the case of censored responses.

Consider two functions $\rho_1$ and $\rho_2$ that satisfy the regularity conditions P1–P4. Moreover, assume that $\rho_2(u) \leq \rho_1(u)$ for all $u$ and that $\sup_u \rho_2(u) = \sup_u \rho_1(u)$. Let $\tilde{\boldsymbol{\beta}}_n$ and $s_n$ be the S-regression and S-scale estimators calculated as in (2.13) and (2.14), respectively. For each $\boldsymbol{\gamma} \in \mathbb{R}^p$ define $R(\boldsymbol{\gamma})$ as

$$R(\boldsymbol{\gamma}) = \frac{1}{n} \sum_{i=1}^n E_{F^*_{n\tilde{\boldsymbol{\beta}}_n}} \left[ \rho_2\left(\frac{u - \boldsymbol{\gamma}' \mathbf{x}_i}{s_n}\right) \Big| \mathbf{w}_i(\tilde{\boldsymbol{\beta}}_n) \right] \tag{2.16}$$

and let $\tilde{\boldsymbol{\gamma}}_n$ be a local minimum of $R(\cdot)$ such that $R(\tilde{\boldsymbol{\gamma}}_n) \leq R(0)$. The MM-estimator $\hat{\boldsymbol{\beta}}_n$ for censored regression is defined by

$$\hat{\boldsymbol{\beta}}_n = \tilde{\boldsymbol{\beta}}_n + \tilde{\boldsymbol{\gamma}}_n. \tag{2.17}$$

The motivation for the definition in (2.17) is as follows. We improve the initial S-estimator $\tilde{\boldsymbol{\beta}}_n$ by fitting an efficient M-estimator to the residuals of $\tilde{\boldsymbol{\beta}}_n$. The resulting M-estimate $\tilde{\boldsymbol{\gamma}}_n$ is the required correction. Expanding the conditional expectations in (2.16) we obtain

$$R(\boldsymbol{\gamma}) = \frac{1}{n} \sum_{i=1}^n \bigg[ \delta_i \rho_2\left(\frac{r_i(\tilde{\boldsymbol{\beta}}_n) - \boldsymbol{\gamma}' \mathbf{x}_i}{s_n}\right) \tag{2.18}$$
$$+ \frac{(1 - \delta_i)}{1 - F^*_{n\tilde{\boldsymbol{\beta}}_n}(r_i(\tilde{\boldsymbol{\beta}}_n))}$$
$$\times \int_{r_i(\tilde{\boldsymbol{\beta}}_n)}^\infty \rho_2\left(\frac{u - \boldsymbol{\gamma}' \mathbf{x}_i}{s_n}\right) dF^*_{n\tilde{\boldsymbol{\beta}}_n}(u) \bigg].$$

For each $i$ such that $\delta_i = 0$ let $M_i = \{j : r_j(\tilde{\boldsymbol{\beta}}_n) > r_i(\tilde{\boldsymbol{\beta}}_n), \delta_j = 1\}$. Then, we have

$$(2.19) \quad \int_{r_i(\tilde{\boldsymbol{\beta}}_n)}^{\infty} \rho_2\left(\frac{u - \boldsymbol{\gamma}' \mathbf{x}_i}{s_n}\right) dF^*_{n\tilde{\boldsymbol{\beta}}_n}(u) = \sum_{j \in M_i} \rho_2\left(\frac{r_j(\tilde{\boldsymbol{\beta}}_n) - \boldsymbol{\gamma}' \mathbf{x}_i}{s_n}\right) \pi_j$$

and $1 - F^*_{n\tilde{\boldsymbol{\beta}}_n}(r_i(\tilde{\boldsymbol{\beta}}_n)) = \sum_{j \in M_i} \pi_j$, where $\pi_j$, $j \in M = \{j : \delta_j = 1\}$, are the probabilities given to the uncensored $r^*_j(\tilde{\boldsymbol{\beta}}_n)$ by the KM estimator $F^*_{n\tilde{\boldsymbol{\beta}}_n}$. For $i, j = 1, \ldots, n$ let

$$(2.20) \quad \pi_{ij} = \begin{cases} \pi_j \Big/ \left( n \sum_{k \in M_i} \pi_k \right), & \text{if } \delta_i = 0 \text{ and } j \in M_i, \\ 1/n, & \text{if } \delta_i = 1 \text{ and } i = j, \\ 0, & \text{otherwise.} \end{cases}$$

Then, from (2.18) and (2.19) we have

$$(2.21) \quad R(\boldsymbol{\gamma}) = \sum_{i=1}^{n} \sum_{j=1}^{n} \rho_2\left(\frac{r_j(\tilde{\boldsymbol{\beta}}_n) - \boldsymbol{\gamma}' \mathbf{x}_i}{s_n}\right) \pi_{ij}.$$

Since the $\pi_{ij}$'s do not depend on $\boldsymbol{\gamma}$, a local minimum of $R(\boldsymbol{\gamma})$ will satisfy

$$\sum_{i=1}^{n} \sum_{j=1}^{n} \rho'_2\left(\frac{r_j(\tilde{\boldsymbol{\beta}}_n) - \boldsymbol{\gamma}' \mathbf{x}_i}{s_n}\right) \mathbf{x}_i \pi_{ij} = \mathbf{0}.$$

Similarly to the uncensored case, this equation can be written as

$$(2.22) \quad \sum_{i=1}^{n} \sum_{j=1}^{n} w_{ij} (r_j(\tilde{\boldsymbol{\beta}}_n) - \mathbf{x}'_i \boldsymbol{\gamma}) \mathbf{x}_i = \mathbf{0},$$

where

$$w_{ij} = \frac{\rho'_2((r_j(\tilde{\boldsymbol{\beta}}_n) - \boldsymbol{\gamma}' \mathbf{x}_i)/s_n)}{((r_j(\tilde{\boldsymbol{\beta}}_n) - \boldsymbol{\gamma}' \mathbf{x}_i)/s_n)} \pi_{ij}.$$

Hence, a local minimum of $R(\boldsymbol{\gamma})$ is the weighted least squares estimator for the points $(r_i(\tilde{\boldsymbol{\beta}}_n), \mathbf{x}_j)$ with weights $w_{ij}$, $i, j = 1, \ldots, n$. Equation (2.22) suggests that an iterative reweighted least squares algorithm can be used to find a local minimum of $R(\boldsymbol{\gamma})$. Furthermore, since we need to find a local minimum such that $R(\boldsymbol{\gamma}) < R(\mathbf{0})$, and reweighted least squares iterations reduce the objective function (see Remark 1 to Lemma 8.3 in Huber [6], page 186) we can start this algorithm at $\boldsymbol{\gamma} = \mathbf{0}$.

2.4. *τ-estimators.* Another way to obtain estimators with high breakdown and high efficiency for normal errors with censored responses, is to extend the class of $\tau$-estimators (Yohai and Zamar [23]). These estimators are based on an efficient scale estimator, called $\tau$-scale.

Let $\rho_1 : \mathbb{R} \to \mathbb{R}^+$ and $\rho_2 : \mathbb{R} \to \mathbb{R}^+$ satisfy conditions P1–P4, and let $b = E_F(\rho_1)$. Moreover, to obtain consistent estimators, we will assume that $\rho_1$ and $\rho_2$ satisfy:

P5. $\rho_i$, $i = 1, 2$, are continuous, and if $0 \leq v < w$ with $\rho_2(w) < \sup_u \rho_2(u)$ then $\rho_2(v) < \rho_2(w)$.

P6. $2\rho_2(u) - \rho_2'(u)u \geq 0$.

Given a sample $u_1, \ldots, u_n$ let $s_n$ be the solution of

$$\frac{1}{n}\sum_{i=1}^{n} \rho_1(u_i/s_n) = b$$

and define the $\tau$-scale as

$$\tau_n^2 = s_n^2 \frac{1}{n}\sum_{i=1}^{n} \rho_2(u_i/s_n).$$

The extension of the $\tau$-estimators for censored data follows the same lines as the one for S-estimators but using a $\tau$-scale instead of an S-scale.

More specifically, let $S_n(\boldsymbol{\beta}, \boldsymbol{\gamma})$ be as in (2.11) and define $\tau_n(\boldsymbol{\beta}, \boldsymbol{\gamma})$ by

$$(2.23) \qquad \tau_n(\boldsymbol{\beta}, \boldsymbol{\gamma})^2 = S_n(\boldsymbol{\beta}, \boldsymbol{\gamma})^2 \frac{1}{n}\sum_{i=1}^{n} E_{F_{n\beta}^*}\left[\rho_2\left(\frac{u - \boldsymbol{\gamma}'\mathbf{x}_i}{S_n(\boldsymbol{\beta}, \boldsymbol{\gamma})}\right)\Big| \mathbf{w}_i(\boldsymbol{\beta})\right].$$

Let

$$(2.24) \qquad \hat{\boldsymbol{\gamma}}_n(\boldsymbol{\beta}) = \arg\min_{\boldsymbol{\gamma} \in \mathbb{R}^p} \tau_n(\boldsymbol{\beta}, \boldsymbol{\gamma})$$

and define the $\tau$-estimator $\hat{\boldsymbol{\beta}}_n$ as in (2.9) or (2.10).

2.5. *Alternative representation.* In this section, we show an alternative way of writing the estimating equations that define our estimators for censored data. This alternative representation is most useful when computing these estimators for non-smooth functions $\rho(u)$ (e.g., the least median of squares—LMS). We will also use this representation in our proofs in the Appendix. This approach also lets us understand better the connection between the estimators defined in the previous sections and their uncensored counterparts.

Let $r_1, \ldots, r_n$ be a random sample from a distribution $F$, and let $c_1, \ldots, c_n$ be unobservable censoring variables independent from the $r_i$'s. Suppose that we observe $r_i^* = \min(r_i, c_i)$ and let $\delta_i = I\{r_i \leq c_i\}$ where $I\{A\}$ is the indicator function of the event $A$. The Kaplan–Meier estimator of $F$ assigns positive weights only to noncensored observations. Furthermore, the self-consistency property of the Kaplan–Meier estimator (Efron [4]) implies that, if $\pi_j$ is the probability assigned to $r_j^*$ for $\delta_j = 1$, then

$$\pi_j = \frac{1}{n} + \sum_{r_j^* > r_i^*, \delta_i = 0} \pi_{ij}, \qquad (2.25)$$

where the $\pi_{ij}$'s are given by (2.20). Observe that $\pi_{ij}$ can be interpreted as the proportion of the mass from the censored $i$th observation that is assigned to the $j$th point. Note that the mass $1/n$ of each censored observation $r_i^*$ is distributed among all the uncensored $r_j^* > r_i^*$ with $\delta_j = 1$ proportionally to $\pi_j$.

Suppose now that $r_i^* = r_i^*(\boldsymbol{\beta})$ for $1 \leq i \leq n$ are residuals for some vector of regression parameters $\boldsymbol{\beta}$, let $\mathbf{x}_i$, $1 \leq i \leq n$, be the corresponding vectors of covariates and call $\pi_{\boldsymbol{\beta},ij}$ the values given by (2.20). The censored residual sample can be written as $\mathbf{z}_1 = (r_1^*(\boldsymbol{\beta}), \delta_1, \mathbf{x}_1')', \ldots, \mathbf{z}_n = (r_n^*(\boldsymbol{\beta}), \delta_n, \mathbf{x}_n')'$. Consider the discrete distribution function $H_{n\boldsymbol{\beta}}^*$ that assigns mass $\pi_{\boldsymbol{\beta},ij}$ to the point $(r_j^*(\boldsymbol{\beta}), \mathbf{x}_i)$. Following the same arguments leading to (2.21) it is easy to show that for any function $g : \mathbb{R} \times \mathbb{R}^p \to \mathbb{R}$ we have

$$\frac{1}{n} \sum_{i=1}^n E_{F_{n\boldsymbol{\beta}}^*}[g(u, \mathbf{x}_i) | \mathbf{z}_i] = \sum_{i=1}^n \sum_{j=1}^n g(r_j^*(\boldsymbol{\beta}), \mathbf{x}_i) \pi_{\boldsymbol{\beta},ij} = E_{H_{n\boldsymbol{\beta}}^*}[g(u, \mathbf{x})]. \qquad (2.26)$$

Then, $C_n(\boldsymbol{\beta}, \boldsymbol{\gamma})$ in (2.7) can be written as

$$C_n(\boldsymbol{\beta}, \boldsymbol{\gamma}) = E_{H_{n,\boldsymbol{\beta}}^*}\left[\rho\left(\frac{u - \boldsymbol{\gamma}' \mathbf{x}}{s_n}\right)\right].$$

This formula simplifies some computations. For example, consider the jump function $\rho$ defined in (2.15) and the solution $s_n$ to $E_{H_{n\boldsymbol{\beta}}^*}[\rho(u/s_n)] = 1/2$. Noting that the marginal distribution of the first coordinate of $H_{n\boldsymbol{\beta}}^*$ is $F_{n\boldsymbol{\beta}}^*$, we have that

$$s_n = \underset{H_{n\boldsymbol{\beta}}^*}{\operatorname{median}}(|u|) = \underset{F_{n\boldsymbol{\beta}}^*}{\operatorname{median}}(|u|),$$

and thus iterative algorithms are not required.

The following theorem shows that $H_{n\boldsymbol{\beta}}^*$ is consistent to the true joint distribution function $H(u, \mathbf{x}) = F(u) G(\mathbf{x})$ when $\boldsymbol{\beta} = \boldsymbol{\beta}_0$. Moreover, Theorem A.1 in the Appendix, shows that if $\boldsymbol{\beta}_n \xrightarrow{P} \boldsymbol{\beta}_0$, then $H_{n\boldsymbol{\beta}_n}^* \xrightarrow{P} H(u, \mathbf{x})$.

THEOREM 1. *Let $(y_i^*, \mathbf{x}_i, \delta_i)$, $i = 1, \ldots, n$, be observations from a censored linear regression model as in Section 2, and assume that the errors and censoring variables satisfy* R1 *and* R2 *on page 7. Let $H_{n\boldsymbol{\beta}}^*$ be defined as above. Then $H_{n\boldsymbol{\beta}_0}^*(u, \mathbf{x}) \to H(u, \mathbf{x})$ a.s.*

**3. Breakdown point.** In general, for a sample $\mathbf{Z}_n$ of size $n$, the finite-sample breakdown point (Donoho and Huber [3]) of an estimator $\mathbf{T}_n = \mathbf{T}_n(\mathbf{Z}_n)$ is defined as

$$\epsilon_n^*(\mathbf{T}_n, \mathbf{Z}_n) = \min_{1 \leq k \leq n} \{k/n : \sup \|\mathbf{T}_n(\mathbf{Z}_{k,n}^*) - \mathbf{T}_n(\mathbf{Z}_n)\| = \infty\},$$

where the supremum is taken over all possible samples $\mathbf{Z}_{k,n}^*$ which are obtained by replacing $k$ observations from $\mathbf{Z}_n$ with arbitrary values and $\|\mathbf{T}\|$ is the $L_2$ norm.

Let $\mathbf{Z}_n = (\mathbf{z}_1, \ldots, \mathbf{z}_n)$ be a sample from a censored linear regression model, where $\mathbf{z}_i = (y_i^*, \mathbf{x}_i, \delta_i)$, $\mathbf{x}_i \in R^p$. Assume that the rank of $\{\mathbf{x}_1, \ldots, \mathbf{x}_n\}$ is $p$ and let $q = \max_{\|\theta\|=1} \#\{i : \theta' \mathbf{x}_i = 0\}$. Let $m$ be the number of censored observations in the sample, $m = \sum_{i=1}^n \delta_i$. The following theorems show that a lower bound for the breakdown point of S-, MM- and $\tau$-regression estimators is

$$(3.1) \qquad \gamma = k_0/n,$$

where

$$(3.2) \qquad k_0 = \min\left(n\left(1 - \frac{b}{a}\right) - q - m, n\frac{b}{a} - m\right),$$

$b$ is the right-hand side of equation (2.11) and $a = \sup_u \rho(u)$.

THEOREM 2 (Breakdown point of S-estimators). *Let $S$ be a scale estimating functional based on a function $\rho$ satisfying* P1–P4. *Let $\hat{\boldsymbol{\beta}}_n$ be the S-estimator defined in Section 2.2, then*

$$(3.3) \qquad \epsilon_n^*(\hat{\boldsymbol{\beta}}_n, \mathbf{Z}) \geq \gamma.$$

THEOREM 3 (Breakdown point of MM-estimators). *Let $\hat{\boldsymbol{\beta}}_n$ be the MM estimator defined in Section 2.3 with functions $\rho_1$ and $\rho_2$ satisfying* P1–P4, *$\rho_2 \leq \rho_1$ and $a = \sup \rho_2 = \sup \rho_1$. Then $\epsilon_n^*(\hat{\boldsymbol{\beta}}_n, \mathbf{Z}) \geq \gamma$.*

The following theorem is proved in Salibian-Barrera and Yohai [17].

THEOREM 4 (Breakdown point of $\tau$-estimators). *Let $\hat{\boldsymbol{\beta}}_n$ be the $\tau$-estimator defined in Section 2.4 with loss functions $\rho_1$ and $\rho_2$ satisfying* P1–P6. *Then $\epsilon_n^*(\hat{\boldsymbol{\beta}}_n, \mathbf{Z}) \geq \gamma$.*

Note that the lower bound in (3.1) is maximized when $b/a = (1 - q/n)/2$. The smallest possible value of $q$ is $p - 1$, and in this case, the sample is said to be in general position (Rousseeuw and Leroy [15]). Using the optimal $b/a$ we have

$$\epsilon_n^*(\hat{\boldsymbol{\beta}}_n, \mathbf{Z}_n) \geq \frac{1}{2}\left(\frac{n - p + 1 - 2m}{n}\right).$$

Note that when $n \to \infty$ the right-hand side converges to $1/2 - \lambda$, where $\lambda$ is the probability of censoring. This is in agreement with our discussion in the Introduction, where we mention that the breakdown point of the median may be as small as $1/2 - \lambda$ when there are censored observations. Although in linear regression models with uncensored response variables it is possible to obtain robust regression estimators with asymptotic breakdown point of 0.5, we believe that the loss in breakdown-point observed in the censored case is due to the use of the Kaplan–Meyer estimator that may convert censored observations into outliers. We conjecture that this loss cannot be to reduced, at least when the estimate is defined using the Kaplan–Meyer estimate.

**4. Asymptotic properties.** The next theorem shows a property related to the consistency of the S-estimator defined in Section 2.2.

THEOREM 5. *Let $\rho$ satisfy regularity conditions* P1–P4. *Let the errors $u$ and covariates $\mathbf{x}$ in the linear model (1.1) have joint distribution function $H_0(u, \mathbf{x}) = F_0(u) G(\mathbf{x})$ such that $F_0(u)$ is symmetric and has a unimodal density, and $G(\boldsymbol{\beta}' \mathbf{x} \neq 0) = t > b/a$ for all $\boldsymbol{\beta} \in \mathbb{R}^p$. Assume that* R1 *and* R2 *on page 7 hold, and let $\boldsymbol{\gamma}_n(\boldsymbol{\beta}_0) = \arg\min_{\boldsymbol{\gamma}} S_n(\boldsymbol{\beta}_0, \boldsymbol{\gamma})$, where $S_n(\boldsymbol{\beta}, \boldsymbol{\gamma})$ is defined in* (2.11). *Then $\boldsymbol{\gamma}_n(\boldsymbol{\beta}_0) \xrightarrow[n \to \infty]{a.s.} \mathbf{0}$.*

The same kind of arguments used in the proof of Theorem 5 can be used to prove similar results for MM-estimators as defined in Section 2.3. Note that a complete proof of consistency would require to show that if $\boldsymbol{\beta} \neq \boldsymbol{\beta}_0$ then $\|\hat{\boldsymbol{\gamma}}_n(\boldsymbol{\beta})\|$ remains asymptotically away from zero. We have not been able to prove this. However, in all our numerical experiments this property seems to hold.

We can nonetheless prove the local consistency and asymptotic normality of the M-estimates defined in Section 2.1. The proof is based on Theorem 5.1 in Ritov [13] where the author studies M-estimates for censored regression which solve (2.6). Unfortunately, showing that there exists a sequence of consistent solutions of this equation seems to be very difficult. However, it can be shown that there exists a sequence $\boldsymbol{\beta}_n$ of approximate solutions to this equation which is $\sqrt{n}$-consistent and asymptotically normal. More precisely, under some regularity conditions Ritov [13] shows that there exists a sequence $\boldsymbol{\beta}_n$ such that

$$(4.1) \qquad \frac{1}{n^{1/2}} \sum_{i=1}^{n} E_{F^*_{n\boldsymbol{\beta}_n}} [\psi(u) | \mathbf{w}_i(\boldsymbol{\beta}_n)] \mathbf{x}_i \xrightarrow{P} \mathbf{0}$$

and such that $\sqrt{n}(\boldsymbol{\beta}_n - \boldsymbol{\beta}_0) \xrightarrow{D} N(0, A_\psi^{-1} B_\psi A_\psi)$ where

$$(4.2) \quad A_\psi = \int E(\mathbf{x}\mathbf{x}' | c - \boldsymbol{\beta}_0' \mathbf{x} \geq u) W_\psi(u) W_{\psi_0}(u) P(c - \boldsymbol{\beta}_0' \mathbf{x} \geq u) \, dF_0(u),$$

where $c$ is the censoring variable,

$$W_\psi(u) = \psi(u) - \frac{\int_u^\infty \psi(t)\,dF_0(t)}{1 - F_0(u)},$$

$\psi_0(u) = -f_0'(u)/f_0(u)$ and

(4.3) $\quad B_\psi = \int E(\mathbf{xx}'|c - \boldsymbol{\beta}_0'\mathbf{x} \geq u) W_\psi^2(u) P(c - \boldsymbol{\beta}_0'\mathbf{x} \geq u)\,dF_0(u).$

The following theorem shows a similar result for the estimates defined by (2.9). To simplify the proofs we will only consider the case where the error scale $\sigma$ is known.

THEOREM 6.  *Assume that:*

1. *$\rho$ satisfies P1, P2 and P3 and P4 and is three times continuously differentiable with bounded derivatives. Moreover, there exists $c_0$ such that $\rho(c_0) = \max_u \rho(u)$ and $P(\min(y,c) - \boldsymbol{\beta}'\mathbf{x} < c_0) < 1$ for all $\boldsymbol{\beta}$ in a neighborhood of $\boldsymbol{\beta}_0$;*
2. *the errors $u_i$ have a symmetric and a strictly unimodal density $f_0$ with finite information for location, that is, $\int_{-\infty}^\infty (\frac{f_0'(u)}{f_0(u_0)})^2 f_0(u_0) < \infty$;*
3. *the vector of explanatory variables $\mathbf{x}$ has compact support; and*
4. *the matrix $A$ defined in (4.2) is nonsingular.*

*Then, there exists a sequence $\boldsymbol{\beta}_n$ such that* (i) $\sqrt{n}\boldsymbol{\gamma}_n(\boldsymbol{\beta}_n) \xrightarrow{P} 0$ *and* (ii) $\sqrt{n}(\boldsymbol{\beta}_n - \boldsymbol{\beta}_0) \xrightarrow{D} N(\mathbf{0}, A_\psi^{-1} B_\psi A_\psi^{-1})$, *where $A_\psi$ and $B_\psi$ are defined in (4.2) and (4.3), respectively.*

Consider a differentiable function $\rho(u)$ satisfying P1–P4, and let $\rho' = \psi$ with $\psi(0) = a_0 > 0$. For $c > 0$ let $\rho_c(u) = (c/a_0)\rho(u/c)$ and $\psi_c(u) = \rho_c'(u) = (1/a_0)\psi(u/c)$. Then the functions $\rho_c$ satisfy P1–P4 and $\lim_{c\to\infty} \psi_c(u) = u = \psi^*(u)$. It is possible to show that $A_{\psi_c} \to A_{\psi^*}$ and $B_{\psi_c} \to B_{\psi^*}$. Therefore, when $c \to \infty$ the relative asymptotic efficiency of the proposed M-estimate with respect to the Buckley and James estimate tends to 1. Choosing $c$ large enough, this relative efficiency can be as close to 1 as desired. For example, this can be obtained using $\rho(u) = \rho_T(u)$ Tukey's bi-square function with derivative

$$\psi_T(u) = u(1 - u^2)^2 I(|u| \leq 1),$$

where $I(|u| \leq 1) = 1$ if $|u| \leq 1$ and 0 otherwise.

**5. Computing algorithm.** Computing the estimators proposed in Section 2 requires solving a highly complex optimization problem. In this section, we present an efficient algorithm to compute the S-estimators defined in Section 2.2.

We will follow a widely used strategy to approximate the solution of complex optimization problems in robust statistics. This approach is based on generating a large number $N$ of candidate vectors $\boldsymbol{\beta}_1, \ldots, \boldsymbol{\beta}_N$. One way to generate these candidates is by drawing subsamples of size $p$ from the data and adjusting them. The estimator is then approximated by the best candidate $\hat{\boldsymbol{\beta}}_n$. The number of candidates $N$ required to obtain a good approximation can be determined explicitly as in the uncensored case (Rousseeuw and Leroy [15]). In other words, if $\boldsymbol{\beta}_1, \ldots, \boldsymbol{\beta}_N$ are the resampling candidates described above, the approximated estimator $\hat{\boldsymbol{\beta}}_n$ satisfies $\hat{\boldsymbol{\beta}}_n = \boldsymbol{\beta}_k$, where

$$\hat{\boldsymbol{\gamma}}_n(\boldsymbol{\beta}_k)' \mathbf{A}_n \hat{\boldsymbol{\gamma}}_n(\boldsymbol{\beta}_k) = \min_{1 \leq j \leq N} \hat{\boldsymbol{\gamma}}_n(\boldsymbol{\beta}_j)' \mathbf{A}_n \hat{\boldsymbol{\gamma}}_n(\boldsymbol{\beta}_j).$$

We now turn our attention to the calculation of $\hat{\boldsymbol{\gamma}}_n(\boldsymbol{\beta}_j)$ for each candidate $\boldsymbol{\beta}_j$. Recall that this requires to solve the minimization problem given by (2.12). For each $\boldsymbol{\beta}_j$ consider a large number of candidates for $\boldsymbol{\gamma}$ and set $\hat{\boldsymbol{\gamma}}_n(\boldsymbol{\beta}_j)$ to be the best of these candidates. Note that for each fixed $\boldsymbol{\beta}_j$ if $\boldsymbol{\beta}_r$ is good approximation to the true $\boldsymbol{\beta}$, then the vector $\boldsymbol{\beta}_r - \boldsymbol{\beta}_j$ is a natural candidates for $\hat{\boldsymbol{\gamma}}_n(\boldsymbol{\beta}_j)$. This observation follows by noting that in this case the residuals $r_i(\boldsymbol{\beta}_j)$ will follow a linear regression model with coefficients $\boldsymbol{\beta} - \boldsymbol{\beta}_j$. Then, we approximate $\hat{\boldsymbol{\gamma}}_n(\boldsymbol{\beta}_j)$ by the vector $\boldsymbol{\beta}_r - \boldsymbol{\beta}_j$ satisfying

$$S_n(\boldsymbol{\beta}_j, \boldsymbol{\beta}_r - \boldsymbol{\beta}_j) = \min_{1 \leq i \leq N} S_n(\boldsymbol{\beta}_j, \boldsymbol{\beta}_i - \boldsymbol{\beta}_j),$$

where for each pair $\boldsymbol{\beta}, \boldsymbol{\gamma} \in \mathbb{R}^p$, $S_n(\boldsymbol{\beta}, \boldsymbol{\gamma})$ is the M-scale estimator defined in (2.11).

Note that, in principle, this algorithm requires finding $N^2$ scales $S_n(\boldsymbol{\beta}_j, \boldsymbol{\beta}_i - \boldsymbol{\beta}_j)$, $i, j = 1, \ldots, n$. However, this is not always necessary. Suppose that we have already computed $\hat{\boldsymbol{\gamma}}_n(\boldsymbol{\beta}_j)$ for $j = 1, \ldots, i$ and let

$$\kappa_i = \min_{1 \leq j \leq i} \hat{\boldsymbol{\gamma}}_n(\boldsymbol{\beta}_j)' \mathbf{A}_n \hat{\boldsymbol{\gamma}}_n(\boldsymbol{\beta}_j),$$

the best value of the objective function obtained so far. We will need to compute $\hat{\boldsymbol{\gamma}}_n(\boldsymbol{\beta}_{i+1})$ only if

$$\hat{\boldsymbol{\gamma}}_n(\boldsymbol{\beta}_{i+1})' \mathbf{A}_n \hat{\boldsymbol{\gamma}}_n(\boldsymbol{\beta}_{i+1}) < \kappa_i.$$

Divide the set of candidates for $\hat{\boldsymbol{\gamma}}_n(\boldsymbol{\beta}_{i+1})$ into two sets: those with $(\boldsymbol{\beta}_k - \boldsymbol{\beta}_{i+1})' \mathbf{A}_n (\boldsymbol{\beta}_k - \boldsymbol{\beta}_{i+1}) \geq \kappa_i$ (call them $\boldsymbol{\gamma}_1, \ldots, \boldsymbol{\gamma}_{N_1}$) and those with $(\boldsymbol{\beta}_k - \boldsymbol{\beta}_{i+1})' \mathbf{A}_n (\boldsymbol{\beta}_k - \boldsymbol{\beta}_{i+1}) < \kappa_i$ (call them $\tilde{\boldsymbol{\gamma}}_1, \ldots, \tilde{\boldsymbol{\gamma}}_{N_2}$). Note that $\|\hat{\boldsymbol{\gamma}}_n(\boldsymbol{\beta}_{i+1})\| < \kappa_i$ only if

$$\min_{1 \leq j \leq N_1} S_n(\boldsymbol{\beta}_{i+1}, \boldsymbol{\gamma}_j) > \min_{1 \leq j \leq N_2} S_n(\boldsymbol{\beta}_{i+1}, \tilde{\boldsymbol{\gamma}}_j).$$

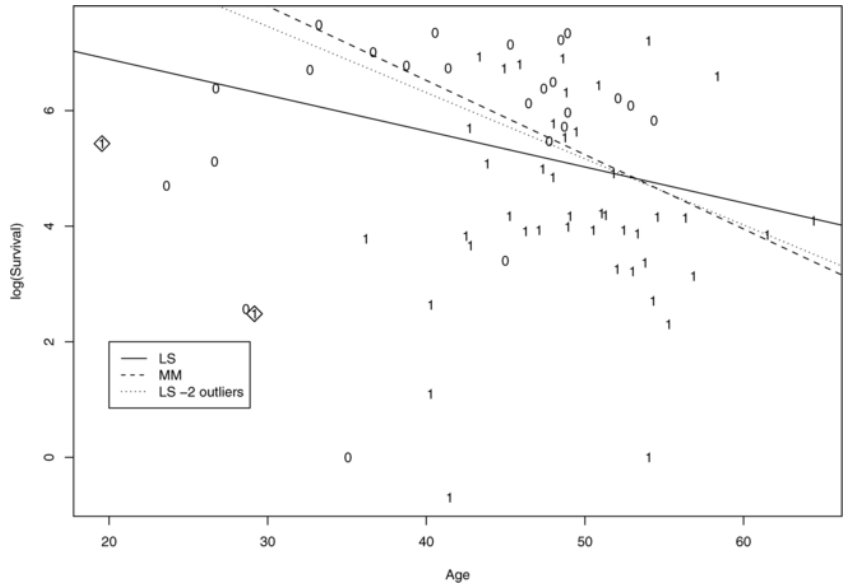

FIG. 2.  *Heart transplant data. "1"s indicate deaths, "0"s indicate censored observations. The least squares estimator seems to be influenced by the two young patients that die early in the study.*

Hence, we first compute $\omega = \min_{1 \leq j \leq N_2} S_n(\boldsymbol{\beta}_{i+1}, \tilde{\boldsymbol{\gamma}}_j)$. Then we compare each $S_n(\boldsymbol{\beta}_{i+1}, \boldsymbol{\gamma}_m)$ for $m = 1, \ldots, N_1$ with $\omega$. If for some $m_0$, we find $S_n(\boldsymbol{\beta}_{i+1}, \boldsymbol{\gamma}_{m_0}) < \omega$ then we stop and set $\kappa_{i+1} = \kappa_i$. Since $\kappa_i \to 0$ we expect $E(N_1)$ to decrease as well. Our Monte Carlo experiments show that there is a substantial gain in speed with this modified algorithm.

**6. Example.** Consider the Heart dataset analyzed in Kalbfleisch and Prentice [8]. These data contain information on heart transplant recipients, including their age and their survival times, which are censored in some cases. In Figure 2, we plot Log (Survival time) versus Age for these patients. We indicate uncensored cases with the symbol "1" and censored ones with "0"s. In the same figure, we also show the fitted lines corresponding to our modified extensions of the LS and MM-estimators. Note that the LS estimator is very much influenced by the early death of two young patients, that can be considered outliers. We used small diamonds around these points to identify them on the plot. We also plot the same LS fit with these two points removed. Note that this line is now close to the robust fit.

**7. Monte Carlo study.** To study the finite-sample properties of these estimators we performed a Monte Carlo study for the simple regression model:

$$y_i = \alpha + \beta x_i + u_i, \qquad i = 1, \ldots, n.$$

TABLE 1
*MSEs without outliers*

| Estimator | S | LMS | LS | MM | GM | L1 |
|---|---|---|---|---|---|---|
| MSE | 0.060 | 0.164 | 0.019 | 0.027 | 0.046 | 0.025 |

We considered 1,000 samples of size $n = 100$, independent normal errors $u_i \sim \mathcal{N}(0,1)$, random covariates $x_i \sim \mathcal{N}(0,1)$ independent from the errors, $\alpha = 0$ and $\beta = 1.5$. We used censoring random variables $c_1, \ldots, c_n$ that were sampled from an independent random variable with distribution $\mathcal{N}(1,1)$. With these choices we have $P(\delta = 0) = 0.32$.

We included the consistent versions under censoring proposed in this paper of the following estimators: the least squares estimator (LS), the least median of squares (LMS), an S-estimator (S) with 50% breakdown point when there is no censoring in the sample, an MM-estimator (MM) with 95% efficiency under normal errors and no censoring, the L1-estimator (L1) [an M-estimator with $\psi(x) = \text{sign}(x)$], and the GM estimator defined by

$$\sum_{i=1}^{n} E_{F^*_{n\beta}}[\psi_1(u - \alpha(\beta))|\mathbf{w}_{\beta i}]\psi_2(x_i - m_x) = \mathbf{0},$$

where $\psi_1(x) = \psi_2(x) = \text{sign}(x)$, $\alpha(\beta) = \text{median}(F^*_{n\beta})$ and $m_x = \text{median}(x_1, \ldots, x_n)$. This is the analogous to the Mood–Brown estimator with breakdown point $1/4$. Both the S- and the MM-estimators used $\rho$ functions in the bisquare family.

The samples were contaminated with 10% of outliers (10 observations). These 10 observations were changed to the points $(x_0, mx_0)$ where $x_0$ was set at 1 and 10 (resulting in low and high leverage outliers resp.), and $m$ ranged between 2 and 5.

In Table 1, we report the MSE for $\beta$ when there are no outliers in the sample. Tables 2 and 3 contain the MSE's for $\beta$ for the cases $x_0 = 1$ and $x_0 = 10$, respectively. From Table 1, we see that, as expected, the most efficient estimator is the LS, followed by the L1 and the MM with efficiencies of 76% and 70%, respectively. For low leverage contaminations (Table 2), the two estimators that perform better, from a maximum MSE point of view, are the L1 and the MM. These two estimators have a similar behavior with a small advantage of the MM. The other estimators are notably worse. Table 3 shows that for high-leverage outliers the MM estimator had the smallest MSE, followed by the S-estimator. Not surprisingly, both the LS and L1 estimators have noticeably worse MSEs than all the other estimators considered here.

Based on these results, we may conclude that the MM-estimators have the best overall performance.

## APPENDIX: PROOFS

### A.1. Consistency of $H^*_{n\beta_0}$.

PROOF OF THEOREM 1. Fix $(a, \mathbf{v}')' \in \mathbb{R}^{p+1}$ and note that $H^*_{n\beta_0}(a, \mathbf{v}) = E_{H^*_{n\beta_0}}[I(u \leq a, \mathbf{x} \leq \mathbf{v})]$ where $I(A)$ denotes the indicator function of the event $A$. Let $u^*_i = y^*_i - \beta'_0 \mathbf{x}_i$ for $i = 1, \ldots, n$. Using (2.26), we have

$$H^*_{n\beta_0}(a, \mathbf{v}) = E_{H^*_{n\beta_0}}[I(u \leq a, \mathbf{x} \leq \mathbf{v})]$$

$$= \frac{1}{n} \sum_{i=1}^{n} \{\delta_i I(u^*_i \leq a, \mathbf{x}_i \leq \mathbf{v})$$

$$+ (1 - \delta_i) E_{F^*_{n,\beta_0}}[I(u_i \leq a, \mathbf{x}_i \leq \mathbf{v}) | u_i > u^*_i]\}.$$

Adding and substracting

$$\sum_{i=1}^{n} (1 - \delta_i) E_H[I(u_i \leq a, \mathbf{x}_i \leq \mathbf{v}) | u_i > u^*_i, \mathbf{x}_i]$$

TABLE 2
MSEs with 10% of outliers at $x_0 = 1$

| | Slopes | | | | | | |
|---|---|---|---|---|---|---|---|
| Estimator | 2 | 2.5 | 3 | 3.5 | 4 | 4.5 | 5 |
| S | 0.10 | 0.27 | 0.38 | 0.30 | 0.20 | 0.13 | 0.10 |
| LMS | 0.14 | 0.30 | 0.54 | 0.69 | 0.79 | 0.76 | 0.78 |
| LS | 0.03 | 0.05 | 0.10 | 0.15 | 0.23 | 0.33 | 0.43 |
| MM | 0.04 | 0.11 | 0.17 | 0.18 | 0.18 | 0.19 | 0.20 |
| GM | 0.09 | 0.25 | 0.40 | 0.52 | 0.62 | 0.71 | 0.78 |
| L1 | 0.07 | 0.16 | 0.20 | 0.21 | 0.21 | 0.21 | 0.21 |

TABLE 3
MSEs with 10% of outliers at $x_0 = 10$

| | Slopes | | | | | | |
|---|---|---|---|---|---|---|---|
| Estimator | 2 | 2.5 | 3 | 3.5 | 4 | 4.5 | 5 |
| S | 0.25 | 0.50 | 0.34 | 0.20 | 0.11 | 0.08 | 0.10 |
| LMS | 0.31 | 0.45 | 0.58 | 0.65 | 0.49 | 0.40 | 0.38 |
| LS | 0.24 | 0.90 | 1.98 | 3.44 | 5.09 | 6.61 | 7.61 |
| MM | 0.23 | 0.45 | 0.30 | 0.17 | 0.08 | 0.06 | 0.07 |
| GM | 0.15 | 0.39 | 0.56 | 0.69 | 0.79 | 0.92 | 1.08 |
| L1 | 0.25 | 0.93 | 2.04 | 3.59 | 5.63 | 8.08 | 11.03 |

we obtain

$$H^*_{n\boldsymbol{\beta}_0}(a, \mathbf{v}) - H(a, \mathbf{v})$$

$$= \frac{1}{n} \sum_{i=1}^{n} [\tilde{g}(u^*_i, \mathbf{x}_i) - H(a, \mathbf{v})]$$

(A.1)

$$+ \frac{1}{n} \sum_{i=1}^{n} [(1 - \delta_i)(E_{F^*_{n,\boldsymbol{\beta}_0}}[g(u_i, \mathbf{x}_i)|u_i > u^*_i]$$

$$- E_H(g(u_i, \mathbf{x}_i)|u_i > u^*_i, \mathbf{x}_i))],$$

where

$$\tilde{g}(u^*_i, \mathbf{x}_i) = \delta_i I(u^*_i \leq a, \mathbf{x}_i \leq \mathbf{v}) + (1 - \delta_i) E_H[I(u_i \leq a, \mathbf{x}_i \leq \mathbf{v})|u_i > u^*_i, \mathbf{x}_i],$$

$H$ denotes the joint distribution of the vector $(\mathbf{x}', u)'$ and $g(u, \mathbf{x}) = I(u \leq a, \mathbf{x} \leq \mathbf{v})$. Note that

$$\tilde{g}(u^*_i, \mathbf{x}_i) = E(I(u_i \leq a, \mathbf{x}_i \leq \mathbf{v})|u^*_i, \mathbf{x}_i, \delta_i)$$

and therefore $E[\tilde{g}(u^*_i, \mathbf{x}_i)] = H(a, \mathbf{v})$.

Since $\tilde{g}$ is bounded, Kolmogorov's law of large numbers yields

$$\frac{1}{n} \sum_{i=1}^{n} [\tilde{g}(u^*_i, \mathbf{x}_i) - H(a, \mathbf{v})] \xrightarrow[n \to \infty]{a.s.} 0.$$

Moreover, note that since $g(u, \mathbf{x}) = I(u \leq a)I(\mathbf{x} \leq \mathbf{v})$ we have

$$E_{F^*_{n,\boldsymbol{\beta}_0}}[g(u_i, \mathbf{x}_i)|u_i > u^*_i] = I(\mathbf{x}_i \leq \mathbf{v}) E_{F^*_{n,\boldsymbol{\beta}_0}}[d(u_i)|u_i > u^*_i],$$

where $d(u) = I(u \leq a)$. Also, because of the independence between $u$ and $\mathbf{x}$ we have $E_H(g(u_i, \mathbf{x}_i)|u_i > u^*_i, \mathbf{x}_i) = I(\mathbf{x}_i \leq \mathbf{v}) E_F(d(u_i)|u_i > u^*_i)$. Hence, the second term in (A.1) equals

$$\frac{1}{n} \sum_{i=1}^{n} (1 - \delta_i) I(\mathbf{x}_i \leq \mathbf{v})(E_{F^*_{n,\boldsymbol{\beta}_0}}[d(u_i)|u_i > u^*_i] - E_F(d(u_i)|u_i > u^*_i)).$$

Thus, we only need to show that

(A.2) $$\sup_{b \in \mathbb{R}} |E_{F^*_{n,\boldsymbol{\beta}_0}}[d(u)|u > b] - E_F[d(u)|u > b]| \xrightarrow[n \to \infty]{a.s.} 0.$$

First, note that we only need to consider the supremum over the set $b \leq a$, since

$$E_{F^*_{n,\boldsymbol{\beta}_0}}[d(u)|u > b] = E_F[d(u)|u > b] = 0 \qquad \text{for } b > a.$$

Next, note that $E_F[d(u)|u > b] = (F(a) - F(b))/(1 - F(b))$. Thus, we need to bound

(A.3)
$$\sup_{b \leq a} \left| \frac{F^*_{n,\boldsymbol{\beta}_0}(a) - F^*_{n,\boldsymbol{\beta}_0}(b)}{1 - F^*_{n,\boldsymbol{\beta}_0}(b)} - \frac{F(a) - F(b)}{1 - F(b)} \right|$$
$$= \sup_{b \leq a} \left| \frac{(F^*_{n,\boldsymbol{\beta}_0}(a) - F(a)) - (F^*_{n,\boldsymbol{\beta}_0}(b) - F(b))}{(1 - F^*_{n,\boldsymbol{\beta}_0}(b))(1 - F(b))} \right.$$
$$\left. + \frac{F(a)(F^*_{n,\boldsymbol{\beta}_0}(b) - F(b)) + F(b)(F(a) - F^*_{n,\boldsymbol{\beta}_0}(a))}{(1 - F^*_{n,\boldsymbol{\beta}_0}(b))(1 - F(b))} \right|$$
$$\leq 4 \frac{\sup_b |F^*_{n,\boldsymbol{\beta}_0}(b) - F(b)|}{(1 - F^*_{n,\boldsymbol{\beta}_0}(a))(1 - F(a))}.$$

Since we are assuming R1 and R2 on page 7, Corollary 1.3 of Stute and Wang [21] implies

$$\lim_{n \to \infty} \sup_b |F^*_{n,\boldsymbol{\beta}_0}(b) - F(b)| = 0 \quad \text{a.s.}$$

This completes the proof. □

THEOREM A.1. *Let $(y^*_i, \mathbf{x}_i, \delta_i)$, $i = 1, \ldots, n$, be observations from a censored linear regression model as in Section 2, and assume that the errors and censoring variables satisfy* R1 *and* R2 *on page 7. Furthermore, assume that $\boldsymbol{\beta}_n \xrightarrow{P} \boldsymbol{\beta}_0$ and let $H^*_{n\boldsymbol{\beta}}$ be defined as above. Then*

$$H^*_{n\boldsymbol{\beta}_n}(u, \mathbf{x}) \xrightarrow{P} H(u, \mathbf{x}).$$

PROOF. The proof follows the same steps as that of the previous theorem replacing $H^*_{n\boldsymbol{\beta}_0}$ by $H^*_{n,\hat{\boldsymbol{\beta}}_n}$. The only difference is that now we need to show that

$$\sup_b |F^*_{n,\hat{\boldsymbol{\beta}}_n}(b) - F(b)| \xrightarrow[n \to \infty]{P} 0.$$

Lemmas 7.1 and 7.2 in Ritov [13] show that

$$\sup_b |F^*_{n,\hat{\boldsymbol{\beta}}_n}(b) - F(b)| \leq O_p(n^{-1/2}) + O(\|\hat{\boldsymbol{\beta}}_n - \boldsymbol{\beta}_0\|) = o_p(1),$$

because $\hat{\boldsymbol{\beta}}_n \xrightarrow{P} \boldsymbol{\beta}_0$. □

**A.2. Breakdown point of the S-estimator.** Define the M-scale estimator $S(F)$ for any arbitrary distribution function $F$ by

$$(A.4) \qquad S(F) = \inf\{s > 0 : E_F[\rho(x/s)] < b\},$$

where $b \geq 0$ and $\rho : \mathbb{R} \to \mathbb{R}_+$ satisfies P1–P4 in Section 2. The following lemma is needed to find the breakdown point of the S-estimators for censored observations. Its proof can be found in Salibian-Barrera and Yohai [17].

LEMMA A.1. *Let $S(F)$ be a scale estimator defined by* (A.4) *where $\rho$ satisfies properties* P1–P4. *Then we have:*

(a) *Given any $K > 0$, and $C > b/a$ there exists $K'$ such that if*

$$(A.5) \qquad P_F\{|x| > K'\} > C,$$

*then $S(F) > K$.*

(b) *Given any $M > 0$ and $C < b/a$, there exist $M'$ such that if*

$$(A.6) \qquad P_F\{|x| > M) < C,$$

*then $S(F) < M'$.*

Given a distribution function $H$ and a Borel set $B$, in the rest of the paper we will denote by $H(B)$ the probability of $B$ under $H$, that is $H(B) = P_H(B)$.

PROOF OF THEOREM 2. Observe that $S_n(\boldsymbol{\beta}, \boldsymbol{\gamma})$ can be defined by

$$(A.7) \qquad E_{H^*_{n,\boldsymbol{\beta}}}(\rho((r - \boldsymbol{\gamma}'\mathbf{x})/S_n(\boldsymbol{\beta}, \boldsymbol{\gamma}))) = b,$$

and $S_n(\boldsymbol{\beta}, \mathbf{0})$ by

$$(A.8) \qquad E_{F^*_{n,\boldsymbol{\beta}}}(\rho(r/S_n(\boldsymbol{\beta}, \mathbf{0}))) = b.$$

Assume that (3.3) is not true. Then there exists a sequence of samples $\mathbf{Z}^{(j)} = (\mathbf{z}_1^{(j)}, \ldots, \mathbf{z}_n^{(j)})$, $1 \leq j < \infty$, $\mathbf{z}_i^{(j)} = (y_i^{*(j)}, \mathbf{x}_i^{(j)}, \delta_i^{(j)})$ such that each $\mathbf{Z}^{(j)}$ differs from $\mathbf{Z}$ in $t$ observations where $t$ satisfies $t < k_0$, and such that if we call $\boldsymbol{\beta}_n^{(j)} = \hat{\boldsymbol{\beta}}_n(\mathbf{Z}^{(j)})$, then

$$(A.9) \qquad \lim_{j \to \infty} \|\boldsymbol{\beta}_n^{(j)}\| = \infty.$$

Let $\boldsymbol{\gamma}_j(\boldsymbol{\beta})$ denote the function $\boldsymbol{\gamma}(\boldsymbol{\beta})$ defined in (2.8) when the sample is $\mathbf{Z}^{(j)}$. We will show that (A.9) is not possible by proving that

$$(A.10) \qquad \lim_{j \to \infty} \|\boldsymbol{\gamma}_j(\boldsymbol{\beta}_n^{(j)})\| = \infty$$

and that

(A.11) $$\sup_j \|\boldsymbol{\gamma}_j(\mathbf{0})\| < \infty.$$

Let us start by proving (A.11). Assume that it is not true. Then without loss of generality we can assume that

(A.12) $$\lim_{j \to \infty} \|\boldsymbol{\gamma}_j(\mathbf{0})\| = \infty$$

and that

(A.13) $$\lim_{j \to \infty} \frac{\boldsymbol{\gamma}_j(\mathbf{0})}{\|\boldsymbol{\gamma}_j(\mathbf{0})\|} = \boldsymbol{\lambda}.$$

We will show that this is not possible by proving that

(A.14) $$\lim_{j \to \infty} S_n^{(j)}(\mathbf{0}, \boldsymbol{\gamma}_j(\mathbf{0})) = \infty$$

and

(A.15) $$\sup_j S_n^{(j)}(\mathbf{0}, \mathbf{0}) < \infty,$$

where $S_n^{(j)}(\boldsymbol{\beta}, \boldsymbol{\gamma})$ denotes the function $S_n(\boldsymbol{\beta}, \boldsymbol{\gamma})$ when the sample is $\mathbf{Z}^{(j)}$.

Let $F_{n,\boldsymbol{\beta},\boldsymbol{\gamma}}^{*(j)}$ denote the distribution of $r - \boldsymbol{\gamma}'\mathbf{x}$ when $(r, \mathbf{x})$ has distribution $H_{n,\boldsymbol{\beta}}^*$ and the sample is $\mathbf{Z}^{(j)}$. Let

(A.16) $$M = \max_{1 \leq i \leq n} |y_i^*| + 1.$$

Then the $y_i^{(j)*}$'s in $\mathbf{Z}^{(j)}$ that are neither contaminated nor censored will have absolute value smaller than $M$. Moreover, $F_{n,\mathbf{0},\mathbf{0}}^{*(j)}$ gives at least mass $1/n$ to each of these points. Therefore, $F_{n,\mathbf{0},\mathbf{0}}^{*(j)}(|y| < M) \geq (n-m-t)/n$. Since $t < k_0$, using (3.2) it follows that $(n-m-t)/n > 1 - b/a$. Thus, from Lemma A.1(b) there exists $M'$ such that $S_n^{(j)}(\mathbf{0}, \mathbf{0}) < M'$ for all $j$, and (A.15) holds.

We now turn our attention to (A.14). Let $\xi_i = |\boldsymbol{\lambda}'\mathbf{x}_i|$, $1 \leq i \leq n$, where $\boldsymbol{\lambda}$ is defined in (A.13), and let

(A.17) $$\xi = \min\{\xi_i : \xi_i > 0\}/2.$$

Then, for all the elements of the original sample, except at most $q$, we have $|\boldsymbol{\lambda}'\mathbf{x}_i| > \xi$. All the contaminated samples $\mathbf{Z}^{(j)}$ have at least $n - q - m - t$ noncensored observations from the original sample $\mathbf{Z}$ such that $|\boldsymbol{\lambda}'\mathbf{x}_i^{(j)}| > \xi$. Then, for $j$ large enough, at least $n - q - m - t$ observations in $\mathbf{Z}^{(j)}$ satisfy

(A.18) $$|y_i^{(j)} - \boldsymbol{\gamma}_j(\mathbf{0})'\mathbf{x}_i| \geq \left| \|\boldsymbol{\gamma}_j(\mathbf{0})\| \left| \left(\frac{\boldsymbol{\gamma}_j(\mathbf{0})}{\|\boldsymbol{\gamma}_j(\mathbf{0})\|}\right)' \mathbf{x}_i \right| - M \right|.$$

Fix $K > 0$ arbitrary and let $K'$ be as in Lemma A.1(a) with $C$ any real number satisfying

(A.19) $$\frac{h_0}{n} > C > \frac{b}{a},$$

where $h_0$ is the smallest integer larger than $nb/a$. Since $t < k_0$, by (3.2) we have $(n - q - m - t)/n > b/a$, and then

(A.20) $$(n - q - m - t)/n > C.$$

Because of (A.12) and (A.13), we can always find $j_0$ large enough so that the right-hand side of (A.18) is larger than $K'$ for all $j > j_0$. Moreover, $F_{n,\mathbf{0},\boldsymbol{\gamma}_j(\mathbf{0})}^{*(j)}$ gives at least mass $1/n$ to those residuals $y_i^{(j)} - \boldsymbol{\gamma}_j(\mathbf{0})'\mathbf{x}_i$. Hence, by (A.20), for $j > j_0$ we have

$$F_{n,\mathbf{0},\boldsymbol{\gamma}_j(\mathbf{0})}^{*(j)}(|y| > K') \geq (n - q - m - t)/n > C.$$

From Lemma A.1(a) it follows that $S_n^{(j)}(\mathbf{0}, \boldsymbol{\gamma}_j(\mathbf{0})) > K$ for all $j > j_0$ and this proves (A.14).

We now prove (A.10). Assume that it is not true. Then we would have

(A.21) $$\sup_j \|\boldsymbol{\gamma}_j(\boldsymbol{\beta}_n^{(j)})\| = L < \infty.$$

To show that this is not possible we will prove that

(A.22) $$\lim_{j \to \infty} S_n(\boldsymbol{\beta}_n^{(j)}, \boldsymbol{\gamma}_j(\boldsymbol{\beta}_n^{(j)})) = \infty$$

and

(A.23) $$\sup_j S_n(\boldsymbol{\beta}_n^{(j)}, -\boldsymbol{\beta}_n^{(j)}) < \infty.$$

To show (A.23) let $M$ be as in (A.16) and observe that there are at least $n - m - t$ observations in $\mathbf{Z}^{(j)}$ with $|y_i^{(j)*}| < M$. It is easy to see that $F_{n,\boldsymbol{\beta}_n,-\boldsymbol{\beta}_n}^{*(j)}$ gives mass at least $1/n$ to these observations, and the proof follows as that of (A.15) above.

We will now prove (A.22). Without loss of generality assume that

(A.24) $$\lim_{j \to \infty} \frac{\boldsymbol{\beta}_n^{(j)}}{\|\boldsymbol{\beta}_n^{(j)}\|} = \boldsymbol{\lambda}.$$

Let $\xi$ be as defined in (A.17). Then for all the elements of the original sample, except at most $q$, we have $|\boldsymbol{\lambda}'\mathbf{x}_i| > \xi$. All the contaminated samples $\mathbf{Z}^{(j)}$ have at least $n - q - m - t$ noncensored observations from the original

sample $\mathbf{Z}$ with $|\boldsymbol{\lambda}'\mathbf{x}_i^{(j)}| > \xi$. Then, for $j$ large enough, at least $n - q - m - t$ observations in $\mathbf{Z}^{(j)}$ satisfy

$$(A.25) \qquad |y_i^{(j)} - \alpha^{(j)\prime}\mathbf{x}_i| \geq \left|\|\beta_n^{(j)}\|\left(\frac{\alpha^{(j)}}{\|\beta_n^{(j)}\|}\right)'\mathbf{x}_i\right| - M\right|,$$

where $\boldsymbol{\alpha}^{(j)} = \boldsymbol{\beta}_n^{(j)} + \boldsymbol{\gamma}_j(\boldsymbol{\beta}_n^{(j)})$. From (A.9), (A.21) and (A.24) it is easy to see that $\lim_{j\to\infty} \boldsymbol{\alpha}^{(j)}/\|\boldsymbol{\beta}_n^{(j)}\| = \boldsymbol{\lambda}$. Observing that $F^{*(j)}_{n,\boldsymbol{\beta}_n^{(j)},\boldsymbol{\gamma}(\boldsymbol{\beta}_n^{(j)})}$ gives at least mass $1/n$ to these $n - m - q - t$ residuals of the form $y_i^{(j)} - \boldsymbol{\alpha}^{(j)\prime}\mathbf{x}_i$, and that the right-hand side of (A.25) can be made arbitrarily large, the rest of the proof follows the same lines as that of (A.14). $\square$

**A.3. Breakdown point of MM-estimators.** The following theorem is needed to find the breakdown point of MM-estimators when the response variable can be censored.

THEOREM A.2. *Let $\mathbf{Z} = (\mathbf{z}_1, \ldots, \mathbf{z}_n)$ with $\mathbf{z}_i = (y_i^*, \mathbf{x}_i, \delta_i)$ and $\mathbf{x}_i \in R^p$ be a sample from a censored linear regression model. Let $\widehat{\beta}_{1n}$ be any regression estimator, and let $\widehat{F}_n^* = F^*_{\widehat{\beta}_{1n},n}$ the KM estimator of the corresponding residual distribution. Let $\rho_1$ and $\rho_2$ two functions satisfying P1–P4, and such that $\rho_2 \leq \rho_1$ and $a = \sup \rho_2 = \sup \rho_1$. Define $s_n = S(\widehat{F}_n^*)$, where $S$ is a M-scale functional based on $\rho_1$ and $0 < b < a$. Let $\widehat{\beta}_{2n}$ be another estimator satisfying*

$$(A.26) \qquad E_{H_n^*}(\rho_2((u + (\widehat{\beta}_{1n} - \widehat{\beta}_{2n})'\mathbf{x})/s_n)) \leq E_{H_n^*}(\rho_2(u/s_n)).$$

*Assume that the rank of $\{\mathbf{x}_1, \ldots, \mathbf{x}_n\}$ is $p$, let $q = \max_{\|\theta\|=1} \#\{i : \theta'\mathbf{x}_i = 0\}$ and $m = \sum_{i=1}^n \delta_i$. Then*

$$(A.27) \quad \epsilon_n^*(\widehat{\boldsymbol{\beta}}_{2n}, \mathbf{Z}) \geq \min(\epsilon_n^*(\widehat{\boldsymbol{\beta}}_{1n}, \mathbf{Z}), (1 - b/a) - (q + m)/n, b/a - m/n).$$

PROOF. Let $\varepsilon_0$ be the right-hand side of (A.27) and assume that the theorem is not true. Then there exists a sequence of samples $\mathbf{Z}^{(j)} = (\mathbf{z}_1^{(j)}, \ldots, \mathbf{z}_n^{(j)})$, $1 \leq j < \infty$, $\mathbf{z}_i^{(j)} = (y_i^{*(j)}, \mathbf{x}_i^{(j)}, \delta_i^{(j)})$ such that each $\mathbf{Z}^{(j)}$ differs from $\mathbf{Z}$ in $t < \varepsilon_0 n$ observations and such that $\lim_{j\to\infty} \|\boldsymbol{\beta}_{2n}^{(j)}\| = \infty$. Since $t < \epsilon_n^*(\widehat{\boldsymbol{\beta}}_{1n}, \mathbf{Z})n$ we have $\sup_j \|\widehat{\boldsymbol{\beta}}_{1n}^{(j)}\| < \infty$. Hence, if we call $\boldsymbol{\gamma}_n^{(j)} = \boldsymbol{\beta}_{1n}(\mathbf{Z}^{(j)}) - \boldsymbol{\beta}_{2n}(\mathbf{Z}^{(j)})$ then

$$(A.28) \qquad \lim_{j\to\infty} \|\boldsymbol{\gamma}_n^{(j)}\| = \infty.$$

Moreover, in all the samples $\mathbf{Z}^{(j)}, 1 \leq j \leq n$, there are at least $n - t - m > (1 - b/a)n$ noncensored observations from the original sample. Since $\sup_j \|\widehat{\boldsymbol{\beta}}_{1n}^{(j)}\| < \infty$ we have that the residuals $r_i^*(\boldsymbol{\beta}_{1n}^{(j)})$ for these $n - t - m$ observations remain bounded uniformly in $j$. Let $\widehat{F}_n^{*(j)}$ be $\widehat{F}_n^*$ when the sample is $\mathbf{Z}^{(j)}$. Then it is clear that $\widehat{F}_n^{*(j)}$ assigns probability at least $1/n$ to these residuals, and hence by Lemma A.1(b) we have $\sup_j S(\widehat{F}_n^{*(j)}) = S^+ < \infty$. Without loss of generality assume that

$$\lim_{j \to \infty} \frac{\boldsymbol{\gamma}_n^{(j)}}{\|\boldsymbol{\gamma}_n^{(j)}\|} = \boldsymbol{\lambda}. \tag{A.29}$$

Let $M = \max_{1 \leq i \leq n} |y_i^*| + 1$, $\delta_i = |\boldsymbol{\lambda}'\mathbf{x}_i|, 1 \leq i \leq n$, and $\delta = \min\{\delta_i > 0\}/2$. Note that all the contaminated samples $\mathbf{Z}^{(j)}$ have at least $n - q - m - t$ non censored observations $\mathbf{z}_i^{(j)} = (y_i^{(j)}, \mathbf{x}_i^{(j)}, \delta_i^{(j)})$ from the original sample $\mathbf{Z}$ which have $|\boldsymbol{\lambda}'\mathbf{x}_i^{(j)}| > \delta$. Then, since for $j$ large enough

$$|y_i^{(j)} - \boldsymbol{\gamma}_n^{(j)\prime} \mathbf{x}_i| \geq \left| \|\boldsymbol{\gamma}_n^{(j)}\| \left|\left(\frac{\boldsymbol{\gamma}_n^{(j)}}{\|\boldsymbol{\gamma}_n^{(j)}\|}\right)' \mathbf{x}_i\right| - M \right|, \tag{A.30}$$

by (A.29) and (A.28), there are at least $n - q - m - t$ observations in $\mathbf{Z}^{(j)}$ such that $|y_i^{(j)} - \boldsymbol{\gamma}_n^{(j)\prime} \mathbf{x}_i| \to \infty$. Since $n_0 = n - q - m - t > nb/a$ we can choose $bn/n_0 < \mu < a$ and let $M = \rho_2^{-1}(\mu)$. There exists a $j_0$ sufficiently large such that for $j \geq j_0$ these $n_0$ observations satisfy

$$|y_i^{(j)} - \boldsymbol{\gamma}_n^{(j)\prime} \mathbf{x}_i|/S^+ > M.$$

Noting that the distribution function $H_n^*$ assigns at least mass $1/n$ to each of these $n_0$ observations, we can conclude that

$$E_{H_n^*}(\rho_2((u + (\widehat{\boldsymbol{\beta}}_{1n}^{(j)} - \widehat{\boldsymbol{\beta}}_{2n}^{(j)})'\mathbf{x})/s_n)) > \frac{n_0}{n}\rho_2(M) > \frac{n_0}{n}\mu > \frac{n_0}{n}\frac{bn}{n_0} = b. \tag{A.31}$$

On the other hand, by the definition of $s_n$ we have

$$E_{H_n^*}(\rho_2(u/s_n)) \leq E_{H_n^*}(\rho_1(u/s_n)) = b. \tag{A.32}$$

Finally, note that (A.31) and (A.32) contradict (A.26). □

PROOF OF THEOREM 3. Follows immediately from Theorem A.2 □

**A.4. Consistency of the S-regression estimator.** Some auxiliary results are needed to prove our main result in this section (Theorem 5). The following lemma is proved as Lemma 7 in Salibian-Barrera and Yohai [17].

LEMMA A.2. *Let $\rho$ satisfy regularity conditions P1–P4. Let $H_n(u, \mathbf{x}) \to F_0(u)G_0(\mathbf{x}) = H_0$ a.s. where $F_0$ is symmetric and has a unimodal density, and $G(\boldsymbol{\beta}'\mathbf{x} \neq 0) \geq t$ for all $\boldsymbol{\beta} \in \mathbb{R}^p$. Then for any $s > 0$ and any $b^* < ta$ there exists $K$ such that*

$$\lim_{n \to \infty} \inf_{\|\boldsymbol{\beta}\| > K} E_{H_n}(\rho((u - \boldsymbol{\beta}'\mathbf{x})/s)) > b^* \qquad a.s.$$

The next lemma is proved as Lemma 9 in Salibian-Barrera and Yohai [17].

LEMMA A.3. *Let $\rho$ satisfy regularity conditions* P1–P4. *Let $H_n(u, \mathbf{x}) \to F_0(u)G_0(\mathbf{x}) = H_0$ a.s. where $F_0$ is symmetric and has a unimodal density, and $G(\boldsymbol{\beta}'\mathbf{x} \neq 0) > b/a$ for all $\boldsymbol{\beta} \in \mathbb{R}^p$. Let $s_0$ be defined by $E_{F_0}(\rho(u)/s_0) = b$. Then given $\varepsilon > 0$ and $K$ there exist $s_1 > s_0$ and $b_1 > b$ such that $\lim_{n \to \infty} \inf_{\varepsilon \leq \|\boldsymbol{\beta}\| \leq K} E_{H_n}(\rho((u - \boldsymbol{\beta}'\mathbf{x})/s_1) > b_1$.*

The next lemma is proved as Lemma 10 in Salibian-Barrera and Yohai [17].

LEMMA A.4. *Let $\rho$ satisfy regularity conditions* P1–P4. *Let $H_n(u, \mathbf{x}) \to H_0(u, \mathbf{x}) = F_0(u)G_0(\mathbf{x})$ a.s. where $F_0$ is symmetric and has a unimodal density, and $G(\beta'\mathbf{x} \neq 0) = t > b/a$ for all $\beta \in \mathbb{R}^p$. Let $s_0$ defined by $E_{F_0}(\rho(u/s_0)) = b$, then if $s_1 > s_0$ we have $\lim_{n \to \infty} E_{H_n}(\rho(u/s_1)) < b$.*

PROOF OF THEOREM 5. Observe that $S_n(\boldsymbol{\beta}_0, \boldsymbol{\gamma})$ is the value $s$ satisfying $E_{H^*_{n,\boldsymbol{\beta}_0}}(\rho((y - \boldsymbol{\gamma}'\mathbf{x})/s)) = b$. We know by Theorem 1 that $H^*_{n,\boldsymbol{\beta}_0}(u, \mathbf{x}) \to H_0(u, \mathbf{x}) = F_0(u)G_0(\mathbf{x})$ a.s. for all $u$ and $\mathbf{x}$. Define $s_0$ by $E_{H_0}(\rho(u/s_0)) = b$. Then using Lemma A.2 with $s = s_0 + 1$, we can find $K$ such that

$$\liminf_{n \to \infty} \inf_{\|\boldsymbol{\gamma}\| > K} S_n(\boldsymbol{\beta}_0, \boldsymbol{\gamma}) \geq s_0 + 1 \qquad \text{a.s.}$$

Let $\varepsilon > 0$ be arbitrary. For this $\varepsilon$ and the $K$ found above, by Lemma A.3, we can find $s_1 > s_0$ such that

$$\liminf_{n \to \infty} \inf_{\varepsilon \leq \|\boldsymbol{\gamma}\| \leq K} S_n(\boldsymbol{\beta}_0, \boldsymbol{\gamma}) \geq s_1 \qquad \text{a.s.}$$

Take $s_2$ such that $s_0 < s_2 < \min(s_0 + 1, s_1)$. By Lemma A.4 we have that $\lim_n S_n(\boldsymbol{\beta}_0, 0) \leq s_2$ a.s. This implies that, with probability 1, there exists $n_0$ such that for all $n \geq n_0$ we have $\|\boldsymbol{\gamma}_n(\boldsymbol{\beta}_0)\| < \varepsilon$. This proves the theorem. $\square$

**A.5. Asymptotic distribution.** Some auxiliary results are needed to prove Theorem 6. The following lemma is proved as Lemma 12 in Salibian-Barrera and Yohai [17].

LEMMA A.5. *Let $H_n(u)$ with $u \in \mathbb{R}^p$ be a sequence of stochastic processes such that, for each $n$ and each element of the underlying probability space where the processes are defined, $H_n(u)$ is a distribution function. Assume that $H_n(u) \xrightarrow{P} H(u)$ for each $u \in \mathbb{R}^p$, where $H(u)$ is a distribution function on $\mathbb{R}^p$. Let $g : \mathbb{R}^p \to \mathbb{R}$ be bounded and continuous, then $E_{H_n}[g(u)] \xrightarrow{P} E_H[g(u)]$.*

The following lemma is proved as Lemma 14 in Salibian-Barrera and Yohai [17].

LEMMA A.6. *Let $\rho$ satisfy regularity conditions* P1–P4. *Let $H_n(u, \mathbf{x}) \xrightarrow{P} F_0(u)G_0(\mathbf{x}) = H_0$ where $F_0$ is symmetric and has a unimodal density, and $G(\boldsymbol{\beta}'\mathbf{x} \neq 0) > t$ for all $\boldsymbol{\beta} \in \mathbb{R}^p$. Assume that $t > E_{F_0}[\rho(u/\sigma)]/a$ where $a = \sup_u \rho(u)$. For all $\varepsilon > 0$ there exists $\delta > 0$ such that*

$$\lim_{n \to \infty} P\left(\inf_{\|\boldsymbol{\alpha}\| > \varepsilon} C_n(\boldsymbol{\beta}_n, \boldsymbol{\alpha}) < E_{F_0}\left(\rho\left(\frac{u}{\sigma}\right)\right) + \delta\right) = 0,$$

*where $C_n(\boldsymbol{\beta}, \boldsymbol{\alpha}) = E_{H^*_{n\boldsymbol{\beta}}}(\rho((u - \boldsymbol{\alpha}'\mathbf{x})/\sigma))$.*

PROOF OF THEOREM 6. Let

$$C_n(\boldsymbol{\beta}, \boldsymbol{\alpha}) = \frac{1}{n}\sum_{i=1}^{n} E_{F^*_{n\boldsymbol{\beta}}}\left[\rho\left(\frac{u - \boldsymbol{\alpha}'\mathbf{x}_i}{\sigma}\right)\bigg|\mathbf{w}_i(\boldsymbol{\beta})\right] = E_{H^*_{n\boldsymbol{\beta}}}\left(\rho\left(\frac{u - \boldsymbol{\alpha}'\mathbf{x}}{\sigma}\right)\right),$$

$$D_n(\boldsymbol{\beta}) = \frac{1}{n\sigma}\sum_{i=1}^{n} E_{F^*_{n\boldsymbol{\beta}}}\left[\psi\left(\frac{u}{\sigma}\right)\mathbf{x}_i\bigg|\mathbf{w}_i(\boldsymbol{\beta})\right] = \frac{1}{\sigma}E_{H^*_{n\boldsymbol{\beta}}}\left(\psi\left(\frac{u}{\sigma}\right)\right),$$

$$L_n(\boldsymbol{\beta}) = \frac{1}{n\sigma^2}\sum_{i=1}^{n} E_{F^*_{n\boldsymbol{\beta}}}\left[\psi'\left(\frac{u}{\sigma}\right)\mathbf{x}_i\mathbf{x}_i'\bigg|\mathbf{w}_i(\boldsymbol{\beta})\right] = \frac{1}{\sigma^2}E_{H^*_{n\boldsymbol{\beta}}}\left(\psi'\left(\frac{u}{\sigma}\right)\mathbf{x}\mathbf{x}'\right).$$

By Theorem 5.1 in Ritov [13], there exists a sequence $\boldsymbol{\beta}_n$ such that

(A.33) $$n^{1/2} D_n(\boldsymbol{\beta}_n) \xrightarrow{p} 0$$

and $n^{1/2}(\boldsymbol{\beta}_n - \boldsymbol{\beta}_0) \xrightarrow{D} \mathcal{N}(\mathbf{0}, A_\psi^{-1} B_\psi A_\psi^{-1})$. Then we only have to prove that $n^{1/2} \gamma_n(\boldsymbol{\beta}_n) \xrightarrow{p} 0$.

Using a second-order Taylor expansion around $\boldsymbol{\alpha} = \mathbf{0}$ we obtain

(A.34) $$C_n(\boldsymbol{\beta}_n, \boldsymbol{\alpha}) = C_n(\boldsymbol{\beta}_n, \mathbf{0}) + D_n'(\boldsymbol{\beta}_n)\boldsymbol{\alpha} + \tfrac{1}{2}\boldsymbol{\alpha}' L_n(\boldsymbol{\beta}_n)\boldsymbol{\alpha} + \|\boldsymbol{\alpha}\|^3 K_n(\boldsymbol{\alpha}),$$

where there exists $\varepsilon_0$ and $K_0$ such that

(A.35) $$p\lim_{n \to \infty} \sup_{\|\alpha\| \leq \varepsilon_0} K_n(\alpha)| \leq K_0.$$

Using Theorem A.1, we have that $H^*_{n\boldsymbol{\beta}_n}(u, \mathbf{x}) \to F_0(u)G_0(\mathbf{x})$ in probability for any $u$ and $\mathbf{x}$, and therefore, by Lemma A.6, we have that for any $\varepsilon > 0$, there exists $\delta > 0$ such that

(A.36) $$\lim_{n \to \infty} P\left(\inf_{\|\boldsymbol{\alpha}\| > \varepsilon} C_n(\boldsymbol{\beta}_n, \boldsymbol{\alpha}) < E_{F_0}\left(\rho\left(\frac{u}{\sigma}\right)\right) + \delta\right) = 0.$$

On the other hand by Lemma A.5

$$C_n(\boldsymbol{\beta}_n, 0) \xrightarrow{P} E_{F_0}\left(\rho\left(\frac{u}{\sigma}\right)\right) = d, \tag{A.37}$$

$$D_n(\boldsymbol{\beta}_n) \xrightarrow{P} E_{H_0}(\psi(u/\sigma))\mathbf{x} = \mathbf{0} \tag{A.38}$$

and

$$L_n(\boldsymbol{\beta}_n) \xrightarrow{P} L_0, \tag{A.39}$$

where

$$L_0 = \frac{1}{\sigma^2} E_{F_0}\left[\psi'\left(\frac{u}{\sigma}\right)\right] E_{G_0}(\mathbf{x}\mathbf{x}'). \tag{A.40}$$

The next step is to prove that $\boldsymbol{\gamma}_n(\boldsymbol{\beta}_n) \xrightarrow{P} 0$. We have

$$\{\|\boldsymbol{\gamma}_n(\boldsymbol{\beta}_n)\| > \varepsilon\} \subset \left\{\inf_{\|\boldsymbol{\alpha}\| > \varepsilon} C_n(\boldsymbol{\beta}_n, \boldsymbol{\alpha}) < d + 2\delta/3\right\} \cup \{C_n(\boldsymbol{\beta}_n, 0) > d + \delta/3\}$$

and therefore (A.36) and (A.37) imply $P\{\|\boldsymbol{\gamma}_n(\boldsymbol{\beta}_n)\| > \varepsilon\}) \to 0$.

Finally, we will prove that $n^{1/2}\|\boldsymbol{\gamma}_n(\boldsymbol{\beta}_n)\| = o_p(1)$. Then if we denote $J_n = \{n^{1/2}\|\boldsymbol{\gamma}_n(\boldsymbol{\beta}_n)\| > \varepsilon\}$, we have to prove that for any $\varepsilon > 0$ we have

$$\lim_{n \to \infty} P(J_n) = 0. \tag{A.41}$$

According to (A.34) we have

$$J_n \subset \left\{\inf_{\varepsilon_0 > \|\alpha\| > \varepsilon n^{-1/2}} [D'_n(\boldsymbol{\beta}_n)\boldsymbol{\alpha} + \tfrac{1}{2}\boldsymbol{\alpha}' L_n(\boldsymbol{\beta}_n)\boldsymbol{\alpha} + \|\boldsymbol{\alpha}\|^3 K_n(\boldsymbol{\alpha})] \leq 0\right\}$$

$$\cup \{\|\boldsymbol{\gamma}_n(\boldsymbol{\beta}_n)\| \geq \varepsilon_0\}.$$

Since $P\{\|\boldsymbol{\gamma}_n(\boldsymbol{\beta}_n)\| > \varepsilon_0\} \to 0$, in order to prove that (A.41) it is enough to show that

$$P\left(\inf_{\varepsilon_0 > \|\alpha\| > \varepsilon n^{-1/2}} \left[\frac{\boldsymbol{\alpha}' D_n(\boldsymbol{\beta}_n)}{\|\boldsymbol{\alpha}\|^2} + \frac{1}{2}\frac{\boldsymbol{\alpha}'}{\|\boldsymbol{\alpha}\|} L_n(\boldsymbol{\beta}_n)\frac{\boldsymbol{\alpha}}{\|\boldsymbol{\alpha}\|} + \|\boldsymbol{\alpha}\| K_n(\boldsymbol{\alpha})\right] > 0\right)$$
$$\to 1 \tag{A.42}$$

and since (A.38), (A.39) and (A.40) hold, it is enough to prove that for all $\varepsilon$

$$p\lim_{n \to \infty} \sup_{\|\alpha\| > \varepsilon n^{-1/2}} \frac{\boldsymbol{\alpha}' D_n(\boldsymbol{\beta}_n)}{\|\boldsymbol{\alpha}\|^2} = 0.$$

This follows from

$$\sup_{\|\alpha\| > \varepsilon n^{-1/2}} \frac{|\boldsymbol{\alpha}' D_n(\boldsymbol{\beta}_n)|}{\|\boldsymbol{\alpha}\|^2} \leq \frac{n^{1/2}}{\varepsilon}\|D_n(\boldsymbol{\beta}_n)\|$$

and (A.33). □

DEPARTMENT OF STATISTICS
UNIVERSITY OF BRITISH COLUMBIA
VANCOUVER, BRITISH COLUMBIA
CANADA V6T 1Z2
E-MAIL: matias@stat.ubc.ca

DEPARTAMENTO DE MATEMÁTICA
UNIVERSIDAD DE BUENOS AIRES
1426 BUENOS AIRES
ARGENTINA
E-MAIL: vyohai@dm.uba.ar